\documentclass[12pt,a4paper,reqno,dvipsnames,color,final]{article}

\usepackage{psfrag,color,graphicx,tikz} 
\usepackage[margin=25mm]{geometry}
\usepackage[utf8]{inputenc}
\usepackage[OT1]{fontenc}

\usepackage{hyperref}

\usepackage{epic,eepic,scrtime,xcolor,mathtools,paralist}
\usepackage[eulergreek]{sansmath}
\usepackage{upgreek}

\usepackage[notcite,notref,color]{showkeys}
\definecolor{refkey}{gray}{.75}   
\definecolor{labelkey}{rgb}{0.3,0,0.6} 

\usepackage{xcolor}

\definecolor{ddorange}{rgb}{1,0.5,0}
\definecolor{ddcyan}{rgb}{0,0.2,1.0}

\usepackage{bef_alex,bm,amsmath,amsfonts,amssymb,mathrsfs}
\newtheorem{hypothesis}[theorem]{Hypothesis}

\def\FG{\bm}
\def\ti{{\times}}

\numberwithin{equation}{section}
\numberwithin{figure}{section}

\makeatletter
\renewcommand*\env@cases[1][1.2]{%
  \let\@ifnextchar\new@ifnextchar
  \left\lbrace
  \def\arraystretch{#1}%
  \array{@{\,}c@{\ }l@{}}%
}
\makeatother

\newcommand{\STEP}[1]{{\underline{\itshape Step #1}}}

\newcommand{\weak}{\rightharpoonup}

\newcommand{\DD}{\calD}  

\newcommand{\dom}{\mathop{\mafo{dom}}}

\newcommand{\ugi}{u^\circ}

\newcommand{\foraa}{\text{for a.a.}}
\newcommand{\Topto}{\stackrel{\Top}{\rightharpoonup}}

\newcommand{\MGS}{(M,\calE,\DD,\psi)}      
\newcommand{\MGStop}{(M,\Top, \calE,\DD,\psi)}   

\newcommand{\BGS}{(X,\calE,\calR)}
\newcommand{\GMGS}{\mathrm{gMGS}}
\newcommand{\GBGS}{\mathrm{gBGS}}
\newcommand{\Spx}{X}

\newcommand{\EL}[2]{\ELname \calE(#1;#2)}
\newcommand{\ELname}{\partial_{\calR}} 

\newcommand{\rsloname}{\mathcal{S}_{\calR}}
\newcommand{\rslo}[1]{\rsloname(#1)}

\newcommand{\wrsloname}{\mathcal{C}_{\calR}}
\newcommand{\wrslo}[2]{\wrsloname(#2;#1)}
\newcommand{\wrslopar}[3]{\mathcal{C}_{\calR_{#1}}(#3;#2)}

\newcommand{\Top}{\mathscr{T}}
\newcommand{\Ls}{\mathop{\mafo{Ls}}}
    
\newcommand{\mlt}[2]{\mathscr{#1}(#2)}

\newcommand{\nameop}{conditioned subdifferential }
\newcommand{\nameslope}{conditioned $\calR$-slope}
\newcommand{\tyos}[3]{\wt{#1}_{#2}^{\,#3}}
\newcommand{\yos}[3]{{#1}_{#2}^{#3}}
\newcommand{\AC}{\mathrm{AC}}
\newcommand{\pairing}[4]{ \sideset{_{#1 }}{_{ #2}}  {\mathop{\langle #3 , #4  \rangle}}}
   
\newcommand{\serifxi}{\upxi}

\newcommand{\llim}[1]{{#1}^{-}}
\newcommand{\rlim}[1]{{#1}^{+}}     
\newcommand{\lrlim}[1]{{#1}^{\pm}}
\newcommand{\Prod}{\mathfrak{P}}

\newcommand{\Jugi}{\boldsymbol{J}}

\newcommand{\RRR}{\color{ddcyan}}

\newcommand{\EEE}{\color{black}}

\renewcommand{\arraystretch}{1.03}
   
\begin{document}

\title{On {De Giorgi}'s lemma for variational interpolants\\
  in metric and Banach spaces%
  \thanks{Research of AM partially supported by DFG via the Berlin Mathematics
    Research Center MATH+ (EXC-2046/1, project ID: 390685689) subproject
    ``DistFell''. Research of RR partially supported by GNAMPA (Indam) and  MIUR via the MIUR-PRIN
    Grant 2020F3NCPX ``Mathematics for industry 4.0 (Math4I4)''.}}

\author{%
 Alexander Mielke\thanks{WIAS Berlin, Anton-Wilhelm-Amo-Str.\,39, 10117 
     Berlin, Germany, alexander.mielke@wias-berlin.de}
 \ and
 Riccarda Rossi\thanks{DIMI, Universit\`a degli studi di Brescia, Italy, 
      riccarda.rossi@unibs.it}}

\date{31 August 2024.  Revised  11 March 2026}

\maketitle

\begin{abstract}
  Variational interpolants are an indispensable tool for the construction of
  gradient-flow solutions via the Minimizing Movement Scheme. The De Giorgi lemma
  provides the associated discrete energy-dissipation inequality. It was
  originally developed for metric gradient systems. Drawing from this theory we
  study the case of generalized gradient systems in Banach spaces, where a
  refined theory allows us to extend the validity of the discrete
  energy-dissipation inequality and to establish it as an equality.  For the
  latter we have to impose the condition of radial differentiability of the
  dissipation potential. Several examples are discussed to show how sharp the
  results are.
\end{abstract} 

\noindent 
\textbf{Keywords:}
Generalized gradient systems, minimizing movement scheme, variational
interpolants, discrete energy-dissipation inequality, radial differentiability.


\bigskip
\medskip

\centerline{\emph{Dedicated to Giuseppe Savar\'e on the occasion of his sixtieth birthday,}}
\smallskip
                \centerline{\emph{with gratitude and friendship}}

\section{Introduction}

The Minimizing Movement Scheme (MMS) was introduced by De Giorgi in \cite{Degi93NPMM}
for constructing solutions for gradient flows in abstract spaces. Since then,
the MMS has developed into a versatile tool for analyzing gradient systems in
Hilbert spaces, Banach space, and metric spaces. In this paper, we address the
 specific tool called ``variational interpolant'', also called ``De Giorgi
interpolant'' that was first introduced in \cite[Lem.\,2.5]{Ambr95MM} and
further developed in \cite{AmGiSa05GFMS}. A
generalization to the Banach spaces was done in \cite[Lem.\,6.1]{MiRoSa13NADN}. 
Variational interpolants  generalize the idea of piecewise affine interpolants
in linear spaces, or geodesic interpolants in geodesic spaces,  in  such a way that they
turn out to be  applicable  in more general situations, namely in general metric
spaces. However, even in the cases of geodesic spaces, including Banach and
Hilbert spaces, they are useful if the energy functional is not geodesically
semi-convex. In general, variational interpolants are no longer
continuous in time and hence, the desired discrete energy-dissipation estimate
is more difficult to obtain. It is exactly this estimate, which is established
in the so-called ``\emph{De Giorgi lemma}".  The purpose of this paper is
twofold: (i) we generalize the validity of the lemma in the Banach setup  and
(ii) we discuss the question why and when the discrete energy-dissipation
estimate is an \emph{equality}. 

To be more precise, we now introduce our approach in more detail by comparing
the theory in metric spaces $(M,\calD)$ and in Banach spaces $(X;\|\cdot\|)$ in
parallel.    It has  to be emphasized that the Banach setup we are going to  address is not
contained in the metric one. Indeed, while in the latter we shall essentially confine the discussion to  the case
in which the dissipation mechanism is encoded by the interplay of the metric $\calD$ with a convex  \emph{scalar} function $\psi:{[0,\infty[}\to {[0,\infty[}$, in the
former we will address \emph{general} dissipation potentials $\calR : X \to
[0,\infty]$, for instance $\calR(v) = \int_\Omega(\frac12 v^2 {+} \frac14 v^4)\dd
x$ for $X=\rmL^4(\Omega)$.   

Following \cite{RMS08} (see also \cite{Miel23IAGS}) we consider a
generalized metric gradient system $\MGS$, subsequently abbreviated by $\GMGS$;
see Definition \ref{def:MGS} for the precise definition.  For a given initial
value $\ugi\in M$ with $\calE(\ugi)<\infty$, it is the aim to construct a curve
$u:{[0,\infty[}\to M$ of maximal slope emanating from $\ugi$, i.e.\ $u$
must satisfy for all $t>0$
\begin{subequations}
\label{eq:I.EDB}
\begin{equation}
  \label{eq:metrCMS}
  \calE(u(t))+ \int_0^t \Big( \psi\big(|u'|(s)\big) +
  \psi^*\big(|\pl\calE|(u(s))\big) \Big) \dd s = \calE(u(0))\quad \text{and}
  \quad u(0)=\ugi,
\end{equation}
where $|u'|\geq 0$ denotes the metric speed of $u$ and $|\pl\calE|(u)\geq 0 $
denotes the metric slope, see \cite{AmGiSa05GFMS},
 and $\psi^*$ is the convex conjugate of $\psi$. 
 The case of general
dissipation functions $\psi:{[0,\infty[}\to {[0,\infty[}$ (lower
semicontinuous, convex, $\psi(0)=0)$, and superlinear) was introduced in
\cite{RMS08}, the choice $\psi(r)=\frac12r^2$ gives the classical notion of
curve of maximal slope of \cite{Ambr95MM}, while $\phi(r)=\frac1p r^p$ leads to
$p$-curves of maximal slopes as in \cite{AmGiSa05GFMS}.

For a generalized Banach-space gradient system $\BGS$, subsequently abbreviated
by $\GBGS$ and precisely defined in Section \ref{ss:3.1}, the aim is to find 
\emph{energy-dissipation balance (EDB) solutions}
$u:{[0,\infty[}\to X$, which are defined via
the following identity
\begin{equation}
  \label{eq:BanEDB}
 \begin{aligned} 
 \calE(u(t))+ \int_0^t\! \Big( \calR\big(u'(s)\big) +
  \calR^*({-} \xi(s) ) \Big) \dd s =  \calE(u(0))& \text{ for all }t>0,
\\
\text{ and } \xi(s) \in \rmD \calE(u(s)) &\text{ for a.a.\ }s\geq 0,
\end{aligned}
\end{equation}
\end{subequations}
where $u'$ is the distributional derivative of $u\in \mathrm{AC}([0,T];X)$, 
  and the convex conjugate $\calR^* : X^* \to [0,\infty)$ is evaluated along
a selection $\xi:{]0,\infty[}\to X^*$  in the multivalued Fr\'echet
subdifferential $\rmD\calE(u)\subset X^*$ of $\calE$, see Section
\ref{ss:3.1}.

With an initial value $\ugi\in M$ and a time step
$\tau>0$ the metric and the Banach-space MMS are defined via $u^0_\tau=\ugi$ and
\begin{subequations}
\label{eq:I.MMS}
\begin{align}
  \label{eq:I.metrMMS}
  u^k_\tau \quad \text{minimizes} \quad M\ni u \ \mapsto \ \tau\;\!
  \psi\big(\tdfrac1\tau \calD(u,u_\tau^{k-1})\big) + \calE(u) \qquad \text{ for all
} k\in \N;
\\
  \label{eq:I.BanMMS}
  u^k_\tau \quad \text{minimizes} \quad X\ni u \ \mapsto \ \tau \,
  \calR\big(\,\tdfrac1\tau (u{-}u_\tau^{k-1})\,\big) + \calE(u) \qquad \text{ for all
} k\in \N.
\end{align}
\end{subequations}
Variational interpolants $\wt u^\tau$ are defined for all $t\in {[0,\infty[}$,
satisfy $\wt u^\tau(k\tau)=u^k_\tau$, and are determined by a variational
condition: for all $k\in \N_0$ and $\sigma\in {]0,\tau[}$, we ask for
\begin{subequations}
\label{eq:I.VarInt}
\begin{align} 
\label{eq:I.metrVarInt}
\wt u^\tau(k\tau{+}\sigma) \quad &\text{minimizes}\quad    M\ni u \ \mapsto \ \sigma
  \psi\big(\tdfrac1\sigma \calD(u,u_\tau^{k-1})\big) + \calE(u);
\\
\label{eq:I.BanVarInt}
\wt u^\tau(k\tau{+}\sigma) \quad &\text{minimizes} \quad \, X\ni u \ \mapsto \ \sigma\,
  \calR\big(\tdfrac1\sigma (u{-}u_\tau^{k-1})\big) + \calE(u).
\end{align}
\end{subequations}
In general, one cannot hope to choose the variational interpolant $t\mapsto \wt
u^\tau(t)$ as a continuous function. However, by classical selection theorems
for measurable multivalued mappings, it is possible to choose a measurable
selection, see Section \ref{ss:3.1}. 

 The  De Giorgi lemma, which was first published in
\cite[Lem.\,2.5]{Ambr95MM}, provides a discrete counterpart to the
energy-dissipation balances in \eqref{eq:metrCMS} and \eqref{eq:BanEDB}, namely
for all $\sigma \in {]0,\tau]}$ we have the so-called \emph{De Giorgi
  estimates}
\begin{subequations}
\label{eq:I.DGLem}
\begin{align}
\label{eq:I.metrDGLem}
&\calE\big(\wt u^\tau(k\tau{+}\sigma )\big) + \sigma\,\psi\big( \tdfrac1\sigma 
 \calD(\wt u^\tau(k\tau{+}\sigma) ,
u^k_\tau)\big) + \int_0^\sigma \psi^* \big( |\pl\calE|(\wt
u^\tau(k\tau{+}\rho) \big) \dd \rho \leq \calE(u^k_\tau); \\
\label{eq:I.BanDGLem}
&\calE(\wt u^\tau(k\tau{+}\sigma )\big) + \sigma\,\calR\big( \tdfrac1\sigma 
 (\wt u^\tau(k\tau{+}\sigma ) {-} u^k_\tau) \big) +  \int_0^\sigma 
 \calR^* \big({-} \xi(k\tau{+}  \rho) \big)  \dd \rho  \leq \calE(u^k_\tau)\\
&\nonumber \hspace*{18em} \text{ for some  } \xi(t)\in \rmD\calE(\wt u^\tau(t))
\text{ for a.a.\ } t>0.
\end{align}
\end{subequations}
For establishing \eqref{eq:I.EDB} via a suitable limit passage, it would be
enough to have \eqref{eq:I.DGLem} for $\sigma=\tau$,  and then adding the
results over all subintervals, but we will see that it is very instructive to
keep $\sigma\in {]0,\tau]}$ on the left-hand side as an independent variable. 

For general differentiable dissipation functions $\psi$,  the De
Giorgi estimate \eqref{eq:I.metrDGLem} was first established in
\cite[Lem.\,4.5]{RMS08}, while $\psi(r)=\frac1pr^p$ is treated in
\cite{AmGiSa05GFMS}. The Banach-space case \eqref{eq:I.BanDGLem} appears first
in \cite[Lem.\,6.1]{MiRoSa13NADN}, but the result therein relies on the
\emph{condition of radial differentiability} of $\calR$, namely
\begin{equation}
  \label{eq:I.RadDiff}
  \forall\, v \in X:\quad \text{the function } {]0,\infty[}\ni \lambda
  \,\mapsto\, \calR(\lambda v) \text{ is differentiable}. 
\end{equation}
This condition is equivalent to the fact that for all
$\xi_1,\xi_2\in \pl\calR(v)$ (with $\partial\calR : X \rightrightarrows X^*$
the convex subdifferential of $\calR$), there holds
$\calR^*(\xi_1)=\calR^*(\xi_2)$, cf.\ Proposition
\ref{prop:equivalence-structural}.

So far, our general overview and introduction shows a complete analogy between
the metric case and the Banach-space setting. Even the condition of radial
differentiability of $\calR$ corresponds to the condition of differentiability
of $\psi$.  However, the methods for establishing the so-called \emph{De
  Giorgi estimates} \eqref{eq:I.metrDGLem} and \eqref{eq:I.BanDGLem} involve
quite different techniques.  In particular, for $\GBGS$ we can exploit the
linear structure of $X$ and thus obtain an Euler-Lagrange equation for the
minimizers $\wt u_\sigma:=\wt u^\tau(k\tau{+}\sigma)$ (keep $k$ fixed,
w.l.o.g.\ $k=0$), namely
\begin{equation}
  \label{eq:I.EulLag}
  0 \in \pl\calR\big(\tdfrac1\sigma (\wt{u}_\sigma{-}u^k_\tau)\big) + \wt \xi_\sigma \quad
  \text{and} \quad  \wt \xi_\sigma \in \rmD\calE(\wt u_\sigma). 
\end{equation}

Indeed, to see a first nontrivial fact, we may assume that $\sigma \mapsto \wt
u_\sigma$ and $\sigma\mapsto \calE(\wt 
u_\sigma )$ are absolutely continuous and such that the chain rule relation 
$\frac\rmd{\rmd\sigma} \calE(\wt u_\sigma )= \langle
\wt\xi_\sigma,\frac\rmd{\rmd\sigma}\wt u_\sigma\rangle$ holds. 
Then, the Euler-Lagrange equation \eqref{eq:I.EulLag} gives the chain rule 
\[
\tdfrac\rmd{\rmd\sigma}\calR\big(\tdfrac1\sigma( \wt 
u_\sigma{-}u^k_\tau)\big) = \big\langle {-} \wt\xi_\sigma,  \frac1\sigma
\tdfrac\rmd{\rmd\sigma}\wt u_\sigma - \frac1{\sigma^2}( \wt u_\sigma{-}u^k_\tau)
\big\rangle.
\]
Thus, differentiating the right-hand side of \eqref{eq:I.BanDGLem} with respect to
$\sigma$ gives 
\begin{align*}
\frac\rmd{\rmd \sigma} \text{\,RHS\eqref{eq:I.BanDGLem}}&= \big\langle
\wt\xi_\sigma , \frac\rmd{\rmd\sigma}\wt u_\sigma\rangle +
\calR\big(\tdfrac1\sigma( \wt u_\sigma{-}u^k_\tau)\big)\\
&\quad  + \big\langle {-}\wt\xi_\sigma , 
\frac\rmd{\rmd\sigma}\wt u_\sigma - \tdfrac1{\sigma}( \wt u_\sigma{-}u^k_\tau)
\big\rangle  +
\calR^*\big({-}\wt\xi_\sigma\big) 
\\&
= \calR\big(\tdfrac1\sigma( \wt u_\sigma{-}u^k_\tau)\big)
+\calR^*\big({-}\wt\xi_\sigma\big)  - \big\langle {-}\wt\xi_\sigma , 
\tdfrac1{\sigma}( \wt u_\sigma{-}u^k_\tau)
\big\rangle   \overset{\text{\eqref{eq:I.EulLag}}}= 0,
\end{align*}
by the Fenchel equivalence $\mu \in \pl\calR(v)\ \Leftrightarrow\
\calR(v)+\calR^*(\mu) = \langle \mu,v\rangle$. 

This observation (which we shall revisit in Section \ref{ss:CR}), motivates our
first main result, see Theorem \ref{th:DGL.RadialDiff}, that  the 
De Giorgi estimate \eqref{eq:I.BanDGLem} is indeed an equality, then called
\emph{De Giorgi identity}. This result is established for all measurable
variational interpolants, under the sole additional assumption of radial
differentiability of $\calR$, cf.\ \eqref{eq:I.RadDiff}. However, for this
version of  the  De Giorgi identity, it is essential to be more
specific with the choice of $\wt\xi_\sigma\in \rmD\calE(\wt u_\sigma)$: one has
to restrict to those $\xi \in \rmD\calE(\wt u_\sigma)$ that minimize
$\calR^*({-}\xi)$ subject to the constraint of satisfying the Euler-Lagrange
equation \eqref{eq:I.EulLag}, see \eqref{MEI-Ban-cond} for the precise
definition. We refer to Example \ref{ex:2-revisited} for a very simple case,
where this restriction is essential for the validity of  the  De Giorgi
estimate as an identity.

For the case of a general $\calR$, dropping radial differentiability, we are
able to establish  the De Giorgi estimate \eqref{eq:I.BanDGLem} if $X$
is a reflexive Banach space.  In fact, along with \eqref{eq:I.BanDGLem}
we shall also obtain a refined estimate, involving a force selection
$\wt \xi_\sigma \in \rmD\calE(\wt u_\sigma)$ that also satisfies the
Euler-Lagrange equation \eqref{eq:I.EulLag}.  In fact, we shall refer to
\eqref{eq:I.BanDGLem} as the \emph{simple De Giorgi estimate}, which will be
enhanced to the \emph{improved De Giorgi estimate} keeping track of
\eqref{eq:I.EulLag}.  Both, the \emph{simple} and the \emph{improved} De Giorgi
estimates will be proved  in our second main result, Theorem
\ref{th:DGL.NON.RD.Simple}.  For this, we use a Yosida-Moreau regularization
$\calR_\eta$ of $\calR$ with an equivalent norm $\| \cdot\|$ such that
$u\mapsto \|u\|^2$ is differentiable. Then, $\calR_\eta$ is differentiable, in
particular also radially differentiable, and for the corresponding $\GBGS$
$(X, \calE, \calR_\eta)$  the  De Giorgi identity holds thanks to
Theorem \ref{th:DGL.RadialDiff}. It can be shown that in the limit passage
$\eta \to 0^+$  the  De Giorgi estimate survives.

While Section \ref{s:3} introduces the definitions and conditions for the case
of gradient systems in Banach space that will then be the focus of Section
\ref{s:4}, we start in Section \ref{se:A.De Giorgi} with the metric case. The
missing Euler-Lagrange equation is replaced by a purely metric identity, not
involving the slope $|\pl\calE|$ but rather the functions $\mathsf d^+$ or
$\mathsf d^-$ defined via
\begin{align*}
 \mathsf d^-_\rho(\ugi)&=\inf\bigset{\DD(\ugi,u)}{u \in
  \mafo{Argmin}\big(\DD(\ugi,\cdot)^2/(2\rho)+\calE(\cdot)\big) },
\\
\mathsf d^+_\rho(\ugi)& =\sup\bigset{\DD(\ugi,u)}{u \in
  \mafo{Argmin}\big(\DD(\ugi,\cdot)^2/(2\rho)+\calE(\cdot)\big) }.
\end{align*}
For $\GMGS$ with differentiable $\psi$, the \emph{metric energy identity} takes
the form
\begin{equation}
  \label{eq:I.metrEnId}
  \calE(\wt u_\sigma) + \sigma\,\psi\big(\tdfrac1\sigma \DD(\ugi,\wt u_\sigma)\big) +
  \int_{\rho=0}^\sigma  \psi^*\Big( \psi'\big(\tdfrac1\rho
  \mathsf{d}^\pm_\rho(\ugi)\big)\Big) \dd \rho = \calE(\ugi) . 
\end{equation}
For $\psi(r)=\frac1pr^p$ the identity was established in
\cite[Thm.\,3.1.4]{AmGiSa05GFMS}, whereas the general case is contained in
\cite[Sec.\,4.2]{Miel23IAGS}.  

Clearly, the metric De Giorgi estimate \eqref{eq:I.metrDGLem} follows easily
from \eqref{eq:I.metrEnId} by inserting the slope inequality
\begin{equation}
  \label{eq:I.SlopeEst}
  |\pl\calE|(\wt u_\sigma) \leq \psi'\big(\tdfrac1\sigma \calD(\wt u_\sigma,
  u^k_\tau)\big) \,\big|\pl({-}\calD(\ugi,\cdot)\big|( \wt u_\sigma) 
\leq \psi'\big(\tdfrac1\sigma \calD(\wt u_\sigma, u^k_\tau)\big),
\end{equation}
see Proposition \ref{pro:MetrSlopeEstim}.  The latter can be seen as a metric
counterpart of the Euler-Lagrange equation.  If one of the two inequalities in
\eqref{eq:I.SlopeEst} is strict, then the  De Giorgi identity is
lost. The last inequality is strict, if $\calD$ is not a length distance, which
means that the metric De Giorgi estimate can only hold in geodesic
spaces. However, even there the first inequality may be strict. In Theorem
\ref{th:MGS.Equality} we show that full continuity and a uniform slope estimate
are sufficient to establish  the  De Giorgi identity in geodesic metric
spaces.

\section{The metric case}
\label{se:A.De Giorgi}

Following \cite{Miel23IAGS}, we specify in the following definition the notion
of metric gradient system we will be working with hereafter.

\begin{definition}
\label{def:MGS}
We call a quadruple $\MGS$ a \emph{generalized  metric gradient system}
(most often abbreviated to $\GMGS$),  if
\begin{compactenum}
\item $(M,\DD)$ is a complete metric space;
\item $ \calE: M \to (-\infty,\infty] $ is a proper (i.e., with non-empty
  domain $\dom(\calE)$) lsc functional;
\item $\psi: \R \to [0,\infty)$ is proper, convex, with $\psi(0)=0$ and
  $ \lim_{r\to \infty} \frac{\psi(r)}{r} = \infty$.
\end{compactenum}
\end{definition}

For later use, we introduce the energy sublevels
\begin{equation}
\label{energy-sublevel}
S_E : = \{ u \in M\, : \ \calE(u) \leq E\}, \qquad E>0.
\end{equation}

\begin{remark}\slshape
Most often, a $ \GMGS$  is in fact individuated by a \emph{quintuple}
$\MGStop $, where, mimicking the setup considered in \cite{AmGiSa05GFMS}, in
addition to the topology induced by the metric $\DD$, a second (Hausdorff)
topology $\Top$ is considered on $M$. Typically, $\Top$ is related to
`coercivity' properties of the energy functional, as it turns out to be the
topology w.r.t.\ which the sublevel sets $S_E$, or the sublevels of a
perturbation of $\calE$, are compact.  Although weaker than the topology
induced by $\DD$, $\Top$ is related to it by the following compatibility
condition (here $\Topto$ denotes convergence with respect to $\Top$):
\[
  (u_n,v_n) \Topto (u,v) \ \Longrightarrow \ \lim_{n\to\infty} \DD(u_n,v_n)
  \geq \DD(u,v).
 \] 
Nonetheless, we have opted for omitting the role of $\Top$ in the discussion of
the metric case in order to avoid overburdening it, on the one hand, and to
highlight the purely metric flavor of the arguments, on the other
hand. Instead, in the Banach setup it will be convenient to encompass the weak
topology in the picture.
\end{remark}

In the setup of a $\GMGS $ $\MGS$, the classical notion of curve of maximal
slope is extended by the following definition (cf.\ \cite[Def.\
4.8]{Miel23IAGS}),
 which brings into play the convex conjugate  $\psi^*: [0,\infty[ \to [0,\infty[$,
$
\psi^*(s) = \sup_{r\in [0,\infty[} (sr {-}\psi(r))
$,  of $\psi$.  
  To simplify the arguments, we fix an arbitrary $T>0$ and
confine the discussion to evolutions on the compact time interval $[0,T]$.
\begin{definition}
\label{def-cms}
Given a generalized  metric gradient system $\MGS$, we say that
$u:[0,T]\to M$ is a \emph{curve of maximal slope} if $u\in \AC ([0,T];M) $ and
it satisfies for every $0 \leq s \leq t \leq T$
\begin{equation}
\label{EDB-cms}
 \calE(u(t)) +\int_s^t \Big( \psi\big( |u'|(r)\big) {+}
\psi^*\big(|\pl\calE|(u(r))\big) \Big) \dd r =  \calE(u(s)) \,. 
\end{equation}
\end{definition}

\begin{remark}\slshape 
\label{rem:AltChainRule}
In \cite[Prob. 2.6]{RMS08} an alternative definition for the above concept was
given, imposing  for the curve $u\in \AC ([0,T];M) $ the \emph{pointwise
  estimate}
\begin{equation}
  \label{def-cms-s2}
  \frac\rmd{\rmd t} \calE(u(t)) \leq - \psi\big( |u'|(t)\big)-
  \psi^*\big(|\pl\calE|(u(t))\big) \qquad \foraa\, t \in (0,T).
\end{equation}
In fact, if the slope $|\pl\calE|$ is a \emph{strong upper gradient} according
to the terminology of \cite{AmGiSa05GFMS} (namely, if a suitable
\emph{chain-rule} inequality holds along $u$), then \eqref{def-cms-s2} is in
fact equivalent to the energy-dissipation balance \eqref{EDB-cms}.
\end{remark}

The Minimizing Movement Scheme for constructing curves of maximal slope
fulfilling the initial condition $u(0)=u_0$, for an assigned initial datum
$u_0 \in M$, then reads as follows: given a time step $\tau>0$, inducing a
(uniform, without loss of generality) partition
$\mathscr{P}_\tau = \{ t_\tau^k \}_{k=1}^{K_\tau}$ of the interval $[0,T]$,
starting from $u_\tau^0:=u_0$ find $(u^{k}_\tau)_{k=1}^{K_\tau}$ such that
\begin{equation}
\label{MMS}
\tag{$\mathrm{MMS}$}
u^{k}_\tau \quad \text{minimizes} \quad M \ni u\mapsto \left( \tau 
\psi\big(\frac1\tau\DD(u^{k-1}_\tau,u)\big) + \calE(u) \right) 
\quad \text{for } k \in \{ 1, \ldots, K_\tau\}\,.
\end{equation}
That is why, from now on we will study the properties of the single-step
minimum problem
\begin{equation}
 \label{eq:single-step-MM}
 \mathop{\mafo{Min}}_{u\in M}\Phi_\sigma(\ugi;u)  \quad \text{with } \  
  \Phi_\sigma(\ugi;u) := \sigma \psi \big(\frac1\sigma\DD(\ugi,u)\big)
  +\calE(u), \quad \sigma >0,
\end{equation}
for a fixed $\ugi \in M$.  We will also use the following notation
\begin{equation}
  \label{resolvent-set}
  \begin{aligned}
   \phi(\ugi;\sigma) =\inf\bigset{\Phi_\sigma(\ugi;u)}{u \in M} \  \text{ and }
   J_\sigma(\ugi)= \mafo{Argmin}\bigset{\Phi_\sigma(\ugi;u)}{u \in M}.
 \end{aligned}
\end{equation}
for the associated value functional and the set of minimizers (which we will
assume non-empty, cf.\ \eqref{non-empy-resolvent} below).  It is also
significant to introduce the following quantities
\[
  \mathsf{d}^-_\sigma(\ugi):= \inf\bigset{ \DD(\ugi,u)}{ u \in J_\sigma(\ugi)}
  \quad \text{and} \quad \mathsf{d}^+_\sigma(\ugi):= \sup \bigset{
    \DD(\ugi,u)}{ u \in J_\sigma(\ugi)}.
\]

Throughout this section, we will work under the following assumptions.

\begin{hypothesis}[Conditions for  generalized  metric gradient systems]
\label{h:M}
We \\ assume that 
\begin{compactitem}
\item[\textbullet] $\psi \in \rmC^1(\R)$ is strictly convex;
\item[\textbullet] $\calE$ is bounded from below by $E_0$, namely $\inf_{u\in
    M} \calE(u)>E_0>0$;
\item[\textbullet] there exists $\sigma_*>0$ such that 
\begin{equation}
\label{non-empy-resolvent}
 J_\sigma(\ugi) 
 \neq \emptyset \text{ for all } \sigma \in (0,\sigma_*) \
 \text{ and all } \ugi \in\dom(\calE).
\end{equation}
\end{compactitem}
\end{hypothesis}

\begin{remark}
\label{rmk:non-empty-resolvent}
\slshape Whenever the generalized metric gradient system  is individuated
by a quintuple $\MGStop $ such that the topology $\Top$ is compatible with
$\DD$, \eqref{non-empy-resolvent} follows from a coercivity property of the
following type: there exists $\sigma_*>0$ such that for all $\sigma \in
(0,\sigma_*)$ and all $\ugi \in M$, for all sequences $(u_n)_n\subset M$ we
have the implication
 \begin{equation}
\label{h:coercivity}
\begin{aligned}
  \sup_n \Phi_{\sigma}(\ugi;u_n) 
  <+\infty \ \Longrightarrow \ (u_n)_n \text{
    admits a $\Top$-converging subsequence}.
\end{aligned}
\end{equation}
Then, the \emph{direct method} yields the existence of minimizers for
\eqref{eq:single-step-MM}.
\end{remark}
 \par
Since this paper is focused on the one-step minimum problem  \eqref{eq:single-step-MM}  and its `starting' point $\ugi$ is fixed,
throughout most of the paper (up to a few exceptions) we will omit to indicate it in the notation for the   functions and sets related to   \eqref{eq:single-step-MM}. Thus, we will  simply write
\[
\Phi_\sigma(u), \quad \phi(\sigma), \quad J_\sigma, \quad   \mathsf{d}^\pm_\sigma\,.
\]

We can in fact enhance \eqref{non-empy-resolvent} by observing that 
\begin{equation}
\label{exist-meas-selec}
\text{there exists a \emph{measurable} selection } (0,\sigma_*) \ni 
\sigma \mapsto \wt{u}_\sigma \in  
J_\sigma
\end{equation}
To show this,  let us consider  the multivalued mapping 
\begin{equation}
\label{Gamma-mapping}
   \Jugi:
   (0,\infty) \rightrightarrows \Spx, \qquad 
    \Jugi(\sigma) := J_\sigma. 
\end{equation}
It is easy to check that $\Jugi$ is upper semicontinuous from $\R$ to
$(M,\calD)$ in that it fulfills for every
$(\sigma_n)_n,\, \sigma \in (0,\infty)$
\begin{subequations}
\label{multivalued-usc}
\begin{align}
&
\label{multivalued-usc-1}
\sigma_n \to \sigma  \ \Longrightarrow \ \Ls_{n\to\infty}
\Jugi(\sigma_n) \subset \Jugi(\sigma)\,, 
\intertext{where  the \emph{Kuratowski upper limit}  (cf., e.g., 
 \cite[Def.\ 4.4.13]{Ambrosio-Tilli})  of the sequence  of closed 
sets $(\Jugi(\sigma_n))_n$ 
  is  defined by }
&
\label{multivalued-usc-2}
u \in  \Ls_{n\to\infty} \Jugi(\sigma_n) \ 
\Longleftrightarrow \ \exists\, (\sigma_{n_k})_k,\,
(u_k)_k \text{ with }
\begin{cases}
\displaystyle 
u_k \in  \Jugi(\sigma_{n_k}) 
 \text{ for all } k \in \N\,,
\\
\lim_{k\to\infty}\calD(u, u_k) =0 \,.
 \end{cases}
\end{align}
\end{subequations}
Since $\Jugi$ is upper semicontinuous, its graph is a Borel subset of
$(0,\infty) \times \Spx$, and the von Neumann-Aumann selection theorem
\cite[Thm.\ 3.22]{Castaing-Valadier77} applies, yielding
\eqref{exist-meas-selec}.

\paragraph{\bf The variational interpolant.}
Thanks to condition \eqref{non-empy-resolvent} the MMS does
admit solutions $(u_\tau^k)_{k=1}^{N_\tau}$.  We are then in a position to
precisely introduce the notion of interpolant of the values
$(u_\tau^k)_{k=1}^{N_\tau}$ we will focus on hereafter.

\begin{definition}[De Giorgi variational interpolant]
  We denote by $\wt u_\tau: [0,T]\to M$ any measurable function obtained by
  setting
\begin{equation}
\begin{aligned}
  & \wt u_\tau(0): = u_\tau^0,
  \\
  & \wt u_\tau(r) \in J_r(u_\tau^{k-1}) =
  \mafo{Argmin}\bigset{\Phi_r(u_\tau^{k-1};u)}{u \in M} \quad \text{if } t =
  t_\tau^{k-1} {+} r.
\end{aligned}
\end{equation}
\end{definition}

The cornerstone of the proof 
of the convergence,
 as $\tau \downarrow 0$,  of (a subsequence of)
the sequence $(\wt u_\tau)_\tau$  to a curve of maximal slope is the
discrete estimate obtained by applying  the  De Giorgi estimate 
\eqref{eq:I.metrDGLem} to the interpolant $\wt u_\tau$. The next section
revolves around the validity of \eqref{eq:I.metrDGLem} as an equality.

\subsection{Metric energy identity and  the De Giorgi estimate:  
statements and examples}
\label{su:2.1}

The following identity was established in \cite[Sec.\,3.1,
eqn.\,(3.1.27)]{AmGiSa05GFMS} for the case $\psi(\delta) = \delta^p/p$ and in
\cite[Thm.\,4.17]{Miel23IAGS} for general differentiable scalar dissipation
potentials $\psi$.

\begin{proposition}[Metric energy identity]
\label{pr:MetrEnergyIdent}
Under Hypothesis \ref{h:M}, any measurable selection
$(0,\sigma_*) \ni \sigma \to \wt{u}_\sigma \in J_\sigma (\ugi) $ fulfills
\begin{equation}
 \label{MEIsec2}
 \calE(\wt u_\sigma) + \sigma\,\psi\big(\frac1\sigma  
  \DD(\ugi,\wt u_\sigma)\big) + \int_0^\sigma \psi^*\big( 
  \psi'\big(\frac1\rho \mathsf{d}^\pm_\rho(\ugi)\big)\big) \dd
 \rho= \calE(\ugi).
\end{equation}
\end{proposition}

There is a straightforward way to relate the metric energy identity to De
Giorgi estimate, and that is throughout the following result in which the
slope at $\wt u_\sigma$ is estimated in terms of the slope of the distance
function $\DD(\ugi,\cdot)$.  In fact, \eqref{slope-of-distance-estimate} below
extends to the setup of a $\GMGS$ $\MGS$, the slope estimate proved in
\cite[Lem.\,3.1.3]{AmGiSa05GFMS} in the quadratic case $\psi(r) = \frac12 r^2$.
For completeness, we also record that a version of
\eqref{slope-of-distance-estimate} was proved in \cite[Lemma 4.4]{RMS08} for
non-differentiable dissipation potentials $\psi$.   Note that this estimate
is the metric counterpart of the Euler-Lagrange equation in the Banach-space
setting, and hence it is less precise for a  general $\GMGS$.

\begin{proposition}[Slope estimate]
  Under Hypothesis \ref{h:M}, 
\label{pro:MetrSlopeEstim} 
\begin{equation}
\label{slope-of-distance-estimate}
\big|\pl\calE \big|(\wt u_\sigma) \leq \psi'\big(\frac1\sigma 
 \DD(\ugi,\wt u_\sigma)\big) \;\big|\pl({-}\DD)(\ugi,\cdot)\big|(\wt u_\sigma) 
  \leq \psi'\big(\frac1\sigma  \DD(\ugi,\wt u_\sigma)\big). 
\end{equation}
\end{proposition}
\begin{proof}
We observe that for  arbitrary
$v\in M$ there holds
\begin{align*}
\calE(\wt u_\sigma)- \calE(v) &= \Phi_\sigma(\wt u_\sigma)- 
\Phi_\sigma(v) - \sigma
\psi\big( \frac1\sigma\DD(\ugi,\wt u_\sigma)\big) +\psi\big( \frac1\sigma\DD( \ugi ,v)\big)
\\
& \stackrel{(1)}\leq \sigma \Big( \psi\big( \frac1\sigma\DD(\ugi,v)\big)
 - \psi\big( \frac1\sigma\DD(\ugi,\wt u_\sigma)\big)\Big)
\\
&\stackrel{(2)}\leq  \psi'\big( \frac1\sigma\DD(\ugi,\wt u_\sigma)\big) \:\Big( {-}\DD(\ugi,\wt u_\sigma) +\DD(\ugi,v)  \Big) . 
\end{align*}
where {\footnotesize (1)} is due to $\wt u_\sigma \in J_\sigma $, whereas
{\footnotesize (2)} follows by the convexity of $\psi$.  Taking the positive
part on both sides (using $\psi'\geq 0$), dividing by $\DD(\wt u_\sigma, v)$,
and taking the limsup for $v\to \wt u_\sigma$ gives the  first estimate in
\eqref{slope-of-distance-estimate}.

By the triangle inequality for the distance $\DD$, for any fixed   $\ugi \in M$
the functions $\DD(\ugi, \cdot) $ and ${-}\DD(\ugi,\cdot)$  have slope less or
equal 1. Therefore, using $|{-}\pl\DD(\ugi,\cdot)|\leq 1$,  the second estimate
in \eqref{slope-of-distance-estimate} follows. 
\end{proof}

We can combine Propositions
\ref{pr:MetrEnergyIdent} and \ref{pro:MetrSlopeEstim} and obtain the De
Giorgi lemma for 
 $\GMGS$  $\MGS$, cf.\ also   \cite[Lemma 4.5]{RMS08}.

\begin{theorem}[De Giorgi lemma for $\GMGS$]
\label{th:DeGiLemmaMGS} 
Under Hypothesis \ref{h:M}, any measurable selection
$(0,\sigma_*) \ni \sigma \mapsto \wt{u}_\sigma \in J_\sigma $ fulfills
\begin{equation}
  \label{eq:MetrDeGiEstimate}
  \calE(\wt u_\sigma) + \sigma\,\psi\big(\frac1\sigma \DD(\ugi,\wt u_\sigma)\big) +
\int_0^\sigma \psi^*\big(\big|\pl\calE \big|(\wt u_\rho) \big) \dd
\rho \leq  \calE(\ugi).
\end{equation}
\end{theorem}

With the aim of improving \eqref{eq:MetrDeGiEstimate} to an equality, taking
into account the slope estimate, it is thus natural to check for cases in which
$|\pl \DD(\ugi,\cdot)| =1$.  This is true for a \emph{geodesic distance}.

\begin{lemma}
\label{l:slope1}
Suppose in addition that the metric space $(M,\DD)$ is a \emph{geodesic space}.
Then, for every $\ugi \in M$ we have $|\pl \DD(\ugi,\cdot)| (u)=1$ for every
$u\neq \ugi$.
\end{lemma}
\begin{proof}
  Clearly, it suffices to show that $|\pl \DD(\ugi,\cdot)| (u)\geq1$. For this,
  let us consider the constant-speed geodesic $\gamma $ connecting $\ugi$ to
  $u$, such that $\DD(\gamma_t,\gamma_s) = (t{-}s) \DD(\ugi,u)$ for all
  $0\leq s \leq t \leq 1$.  Then,
  \begin{equation}
  \label{ok4length}
\begin{aligned}
  |\pl \DD(\ugi,\cdot)| (u) & = \limsup_{v\to u} \frac{(
    \DD(\ugi,u){-}\DD(\ugi,v))^+}{\DD(u,v)}
  \\
  & \geq \limsup_{t \to 0^+} \frac{(
    \DD(\ugi,u){-}\DD(\ugi,\gamma_t))^+}{\DD(u,\gamma_t)} = \lim_{t \to
    0^+}\frac{\DD(\ugi,u) - t \DD(\ugi,u)}{(1{-}t) \DD(\ugi,u)} =1\,.
 \end{aligned}
 \end{equation}
\end{proof}  
 \begin{remark}
\label{rmk:ext2length}
\sl
We highlight that, for the validity of the above result it is sufficient for $(M,\DD)$ to be just a \emph{length space} \cite[Chap.\ 2]{BuragoBuragoIvanov}, namely
with the property that the distance $\DD(\ugi,u)$ between any two points $\ugi, \, u \in M$ can be realized as the infimum of the lengths of the curves joining them. 
Then, for the lower estimate in \eqref{ok4length}
the role of the constant-speed geodesic $\gamma $ connecting $\ugi$ to
  $u$ can be played by a curve with length arbitrarily close to 
  $\DD(\ugi,u)$. 
\end{remark}

In contrast, for the non-geodesic distance
\begin{equation}
\label{nongeodes}
\DD(u,w):=\min\{|u{-}w|_2, R\}  \quad \text{on $\R^n$}
\end{equation} 
with $R>0$ a given constant (indeed, note that $(\R^n, \DD)$ is not a length space, either),   we have $|\pl({\pm}\DD)(\ugi, \cdot)|(u)=0$
whenever $|u{-}\ugi|>R$.

We next provide two examples showing that, without a geodesic metric and
without a continuous slope, we cannot expect \eqref{eq:MetrDeGiEstimate} to
hold as an equality.

\begin{example}[Failure of equality in \eqref{eq:MetrDeGiEstimate}  if  $(M,\DD)$ not geodesic]\slshape
\label{ex:NotGeodesic}
We consider the \emph{quadratic} metric gradient system $(M,\calE,\DD)$ with
\[
M=\R, \quad \calE(u)=\frac12 \,u^2, \quad \DD(u,w)=\min\{|w{-}u|, 1\}, \quad
\psi(\delta) = \delta^2/2.
\]
Starting with $\ugi>1$ and setting $\sigma_*=\big((\ugi)^2{-}1\big)^{-1/2}$ we obtain 
\[
 J_\sigma= \mafo{Argmin}\bigset{ 
  \frac1\sigma \DD(\ugi,u)^2 + \calE(u)}{ u\in \R} = 
\begin{cases}\frac{\ugi}{1+\sigma} & \text{for } \sigma < \sigma_*, \\
\big\{ \frac{\ugi}{1+\sigma} , 0\big\}& \text{for } \sigma=\sigma_*,\\
0 & \text{for } \sigma > \sigma_*. 
\end{cases}
\]
We can calculate all terms in \eqref{eq:MetrDeGiEstimate} and find equality as
long as $\sigma \leq \sigma_*$; but strict inequality holds for $\sigma>\sigma_*$. 

Note that $|\pl \calE|(\wt u_\sigma) =|\pl \calE|(0)=0 $ for $\sigma>\sigma_*$
but
$\psi'\big(\frac1\sigma \DD(\ugi,\wt u_\sigma) \big) = \frac1\sigma
\DD(\ugi,0)= 1/\sigma \gneqq 0$.   Hence, the first  estimate in
Proposition \ref{pro:MetrSlopeEstim} holds as an equality,  but the second
estimate  is strict because of  $\big|\pl({-}\DD)(\ugi,\cdot)\big|(0)=0
 <1$.
\end{example}

Our second counterexample involves a \emph{discontinuous} 
slope functional.

\begin{example}[Failure of equality in \eqref{eq:MetrDeGiEstimate}  if the   slope of $\calE$ is not continuous]\slshape
\label{ex:SlopeNotContin}
We consider the  (again, quadratic) metric gradient system 
\[
M=\R,\quad \calE(u)=\max\{u,0\},\quad \DD(u,w)=|u{-}w|,\quad
\psi(\delta)=\delta^2/2.
\]
Starting at $\ugi=1$ we find the unique variational interpolant $\wt u_\sigma=
\max\{ 1{-}\sigma,0\}$. The curve 
$\sigma \mapsto \wt u_\sigma$
is absolutely continuous but the slope along
the curve is discontinuous, namely
\[
|\pl\calE|(\wt u_\sigma) = 1 \ \text{ for }\sigma \in {[0,1[} 
\quad \text{ and }\quad 
|\pl\calE|(\wt u_\sigma) = 0 \ \text{ for }\sigma \geq 1. 
\]
This time, the  first  estimate in Proposition \ref{pro:MetrSlopeEstim}
is strict for $\sigma\geq 1$, and hence \eqref{eq:MetrDeGiEstimate} is also a
strict inequality for $\sigma>1$.
\end{example}

\subsection{Equality in  the   De Giorgi estimate}
\label{suu:Equality}

The discussion in Section \ref{su:2.1} has highlighted the link between two 
properties (one related to the geometry of the space, the other to the driving
energy  and its slope),  and equality in  the De Giorgi Lemma.  With Theorem
\ref{th:MGS.Equality} we now prove that the joint validity of such properties
is a sufficient, albeit rather strong, condition for a $\GMGS $ to
guarantee that all measurable variational interpolants satisfy estimate
\eqref{eq:MetrDeGiEstimate} with equality.  The idea is that, if $|\pl\calE|$ is upper semicontinuous, points 
where the slope is `too big' are avoided by minimizers.

\begin{theorem}
\label{th:MGS.Equality}
Consider a $\GMGS$ $(M,\calE,\DD,\psi)$ satisfying Hypothesis \ref{h:M} such
that, additionally,
\begin{enumerate}\itemsep-0.4em
 \item $(M,\DD)$ is a geodesic space,  and 
 \item $\calE$ is continuous on $M$,
 \item $|\pl\calE|$ is upper semicontinuous on $M$.
\end{enumerate}
Then, for any measurable selection
$(0,\sigma_*) \ni \sigma \mapsto \wt{u}_\sigma \in J_\sigma  $ the
relation \eqref{eq:MetrDeGiEstimate} holds with equality.
\end{theorem}
\begin{proof}
It suffices to show that the slope estimate from Proposition
\ref{pro:MetrSlopeEstim} holds with equality and combine this with Lemma
\ref{l:slope1}.  Since the upper estimate
$|\pl\calE|(\wt u_\sigma) \leq \psi'\big( \frac1\tau \DD(\ugi,\wt
u_\sigma)\big)$ is clear, it remains to show
$|\pl\calE|(\wt u_\sigma) \geq \psi'\big(\frac1\sigma \DD( \ugi,\wt
u_\sigma)\big)$.

We consider the geodesic curve $[0,1]\ni\theta\mapsto \gamma_\theta$ with $\gamma_0=\ugi$ and
$\gamma_1=\wt u_\sigma$. Since $\gamma_1=\wt u_\sigma $ is a global minimizer
for $\Phi_\sigma(\cdot)$, we can compare
with $u=\gamma_\theta$ and obtain  
\begin{align*}
\calE(\gamma_\theta)-\calE(\gamma_1)& = \calE(\gamma_\theta)-\calE(\wt u_\sigma)\geq \sigma\,\psi\big(\frac1\sigma \DD(\ugi,\wt
u_\sigma)\big) - \sigma\,\psi\big(\frac1\sigma \DD(\ugi,\gamma_\theta)\big) \\
&\geq \psi'\big( \frac1\sigma \DD(\ugi,\gamma_\theta)\big)\,\big(
\DD(\gamma_0,\gamma_1)-\DD(\gamma_0,\gamma_\theta) \big) 
=  \psi'\big( \frac1\sigma \DD(\ugi,\gamma_\theta)\big)\, \DD(\gamma_\theta,
 \gamma_1), 
\end{align*}
where we used the convexity of $\psi$ and the fact that $(\gamma_\theta)_\theta$ is
a geodesic.  Thus, we have 
\[
S_\theta:= \frac{\calE(\gamma_\theta){-} \calE(\gamma_1)}{\DD(\gamma_\theta,
 \gamma_1)} \geq  \psi'\big( \frac1\sigma \DD(\ugi,\gamma_\theta)\big) \text{ \
 for all } \theta \in [0,1).
\]
Now fix $\theta$ and define the function $h_\theta(s)= \calE(\gamma_s)-
S_\theta \calD(\gamma_s,\gamma_1)$ for $s \in [\theta,1]$. Since $\calE$ is continuous, $h_\theta$ is
continuous, too, 
 and satisfies
$h_\theta(\theta)=h_\theta(1)=\calE(\gamma_1)$. Hence, there exists a maximizer
$t_\theta \in (\theta,1)$ of $h_\theta$. For this $t_\theta$ and $s \in
(t_\theta,1]$ we have, by construction, 
\[
 \calE(\gamma_{t_\theta})- \calE(\gamma_s) \geq S_\theta
   (\calD(\gamma_{t_\theta},\gamma_1)){-}\calD(\gamma_s,\gamma_1)\big) =
   S_\theta \calD(\gamma_{t_\theta}, \gamma_s) .  
\]  
Dividing by $\calD(\gamma_{t_\theta}, \gamma_s)$ and taking the limit
$s\searrow t_\theta$ we find
\[
|\pl\calE|(\gamma_{t_\theta}) \geq S_\theta \geq
  \psi' \big(\frac1\sigma\calD(\ugi,\gamma_\theta)\big) .
\]
In the limit $\theta\nearrow 1$ we have $\gamma_{t_\theta}\to
\gamma_1=\wt u_\sigma $, and the upper semicontinuity of $|\pl\calE| $ yields 
\[
|\pl\calE|(\wt u_\sigma)=|\pl\calE|(\gamma_1) \geq \limsup_{\theta\nearrow 1} 
|\pl\calE|(\gamma_{t_\theta}) \geq \lim_{\theta\nearrow 1}  \psi'
\big(\frac1\sigma\calD(\ugi,\gamma_\theta)\big) 
 = \psi'\big(\frac1\sigma\calD(\ugi,\wt u_\sigma )\big),
\]
which is the desired opposite inequality to \eqref{MEIsec2}, and the assertion
follows.   
\end{proof}

\section{Banach case: examples and  preliminary results}
\label{s:3}

\subsection{Setup}
\label{ss:3.1}

We consider a generalized Banach-space gradient system ($\GBGS$, for short)
 $(\Spx, \calE, \calR)$ where the state space
\begin{equation}
\label{space}
\tag{$\mathrm{X}$}
\text{$\Spx$ is a reflexive Banach space.} 
\end{equation} 
 We will thus study the minimum problem 
\begin{equation}
\label{min-Banach}
\mathop{\mathrm{Min}}\limits_{u\in \Spx}\Phi_\sigma(u)  \qquad \text{with } \
\Phi_\sigma(u) := \sigma \calR\Big( \frac1\sigma (u{-}\ugi) \Big) 
+\calE(u), \quad \sigma >0,
\end{equation}
for a fixed  $\ugi \in \dom(\calE) \subset  X $.  

We now collect our working assumptions on $\calE$ and $\calR$, which partly
mirror those collected in Definition \ref{def-cms} and Hypothesis \ref{h:M}. At
the same time, they clearly reflect the underlying Banach setup, involving the
weak topology on $\Spx$, in addition to the norm topology, in conditions
\eqref{HypEneBan} below (cf.\ also Remark \ref{rmk:non-empty-resolvent}). As
already mentioned in the Introduction, in the Banach setting we will allow for
\emph{nonsmooth} energies, and thus work with the \emph{Fr\'echet
  subdifferential} $\pl \calE$ of $\calE$ in place of its G\^ateau derivative
$\rmD \calE$. We recall that the multivalued operator
$\pl \calE: \Spx \rightrightarrows \Spx^* $ is defined at
$u \in \dom(\calE) $ by
\begin{equation}
\label{FrSubdiff}
\xi 
\in \pl \calE(u) 
\quad \text{if and only if} \quad  \calE(w) -\calE(u) \geq \langle
{\xi},{w{-}u} \rangle+ \mathrm{o}(\| w{-}u\|_\Spx) \text{ as } w \to u \,. 
\end{equation}
Then, in \eqref{HypEneBan-2} we ask for closedness of the graph of
$\partial \calE$, w.r.t.\ the weak topology of $\Spx {\times} \Spx^*$, along
sequences with bounded energy.

\begin{hypothesis}[Conditions for  $\GBGS$]
\label{h:X}
We assume that
\begin{compactitem}
\item The dissipation potential
  $\calR: \Spx \to [0,\infty]$
  has a proper domain $\dom(\calR)$ \emph{open} in $\Spx$,   it 
  is lower semicontinuous, convex, and fulfills the following conditions:
  \begin{equation}
    \label{Diss}
    \begin{aligned}
     \calR(0) =0, \qquad  \lim_{\|v \|\to \infty} \frac{\calR(v)}{\|v \|} = \infty. 
     \end{aligned}
  \end{equation} 

\item The energy functional $\calE: \Spx \to (-\infty,\infty]$ is proper,
  bounded from below, and weakly-sequentially lower semicontinuous, i.e.\ for
  all $(u_n)_n \subset \Spx$
  \begin{subequations}
    \label{HypEneBan}
    \begin{align}
      \label{HypEneBan-1}
      u_n  \rightharpoonup u\text{ in } \Spx \  \Longrightarrow  \ 
      \liminf_{n\to \infty} \calE(u_n) \geq \calE(u)\,,
    \end{align}
    and $\partial \calE: \Spx \rightrightarrows \Spx^*$ is closed on energy
    sublevels (cf.\ notation \eqref{energy-sublevel}) , i.e.
    \begin{align}
      \label{HypEneBan-2}
      \forall\, E>0: \ \left\{
        \begin{array}{cr}
          (u_n,\xi_n) \rightharpoonup (u,\xi) &  \text{in } \Spx \times \Spx^*,
          \\
          u_n \in S_E , \ \xi_n \in \partial\calE(u_n)  & \text{for all } n \in \N
        \end{array}
      \right\} \  \Longrightarrow  \ \xi \in \partial\calE(u)\,.
    \end{align}
  \end{subequations}
\end{compactitem}
\end{hypothesis}

 We pin down some elementary, but important, consequences of our
assumptions on $\calR$, for its convex conjugate
\[
  \calR^* : \Spx^* \to [0,\infty],
 \qquad \calR^*(\xi) = \sup_{v \in \Spx} (\langle \xi, v \rangle {-}\calR(v)), 
\]
which is, by construction, lower semicontinuous and convex on $\Spx^*$.

\begin{lemma}
\label{le:calR*}
The functional $\calR^*: \Spx^* \to [0,\infty]$ has domain
$\dom(\calR^*) = \Spx^*$, it fulfills $\calR^*(0)=0$, and the
\emph{coercivity} estimate
 \begin{equation}
 \label{coercivityRstar}
 \exists\, r, K_r>0 \ \forall\, \xi \in \Spx^* \, : \quad \calR^*(\xi) \geq r
 \|\xi \|_* -K_r\,. 
 \end{equation}
  Furthermore, $\calR$ is locally Lipschitz on $\dom(\calR)$. 
\end{lemma} 
\begin{proof}
Since $\calR$ has superlinear growth at infinity, $\calR^*$ has full domain 
(cf., e.g., \cite[Prop.\,2.11,\,2.14]{Brez73}). It can be immediately checked
that $\calR^*(0)=0$.  To show \eqref{coercivityRstar}, we observe that since
$0 $ belongs to the open set $\dom(\calR)$,
 $\calR$ is continuous in $0$ \cite[Chap.\,I, Cor.\,2.5]{ET99} and, thus, 
 there exists $r>0$ such
that $\overline{B_r(0)} \subset \dom(\calR)$ and $\calR$ is bounded
from above on $\overline{B_r(0)}$.
Therefore, 
\[
  \calR^*(\xi) \geq \langle \xi, \ell r w \rangle - \calR(\ell r w) \qquad
  \text{for all } w \in \overline{B_1(0)}, \ \ell \in [0,1]\,.
\]
Taking into account that
$ \calR(\ell r w) \leq \ell \calR(r w)\leq \calR(r w)$, we thus conclude
\[
  r \|\xi\|_* = r \sup_{w\in \overline{B_1(0)}} |\langle \xi, w \rangle | \leq
  \calR^*(\xi) + K \qquad \text{with } K_r = \sup_{v \in \overline{B_r(0)}}
  \calR(v)<\infty\,.
\]
The last assertion follows from   \cite[Chap.\,I, Cor.\ 2.4]{ET99}, which states that $\calR$ is locally Lipschitz on  (the open)
$\dom(\calR) $ if and only if there exists a non-empty convex set over which $\calR$ is bounded from above which, in our case, is the aforementioned ball 
$B_r(0) $. 
\end{proof}

\begin{remark}\slshape 
\label{rmk:threshold}
For most part of the subsequent analysis, the open-domain condition for $\calR$
will be sufficient.  One the one hand it  implies the coercivity
property \eqref{coercivityRstar} for $\calR^*$, and on the other hand it
implies local Lipschitz continuity of $\calR$ in all points of the domain. 
Only  for specific results (see Lemma \ref{l:characteriz-geom2} and Proposition \ref{pr:SuffC.CR.lsc} ahead),  
 we will have to
require that $\calR^*$ has superlinear growth at infinity,
 which in particular implies that $\calR$ has
full domain.

It remains an open problem, though, to deal with the case in which
$\dom(\calR)$ is a general proper subset of $X$, not necessarily open,
which occurs, for instance, when the gradient system $(X,\calE,\calR)$ models
the unidirectional evolution of inelastic processes in solid mechanics like
damage.
\end{remark}

\begin{remark}\slshape
\label{rmk:state-depend-diss}
It is often significant to consider dissipation potentials that also depend on
the state variable, i.e.\ $ \calR= \calR(u,v)$, but this generalization would
be irrelevant for the study of the properties of the single-step minimum
problem \eqref{min-Banach} for a fixed  $\ugi \in X $.  Indeed, in the
state-dependent case the dissipation term would be simply replaced by
$ \calR\left( \ugi, \frac1\sigma (u{-}\ugi)\right) $.
\end{remark}

We emphasize that the closedness condition \eqref{HypEneBan-2} assumes only
\emph{weak} convergence in $\Spx$ on sequences $(u_n)_n$. However, the
additional assumptions $(u_n )_n \subset S_E $ and the existence of a
\emph{bounded} sequence $(\xi_n)_n$ such that $ \xi_n \in \partial\calE(u_n) $
for all $n \in \N$, often grants extra compactness properties to the sequence
$(u_n )_n $. 

 Nonetheless, it would also be possible to bypass \eqref{HypEneBan-2}  by working with the so-called \emph{limiting subdifferential} of $\calE$, defined at a given $u \in  \dom(\calE) $ as the set of all $\xi$  that are weak limits of sequences $(\xi_n)_n$, with $\xi_n \in \partial\calE(u_n)$ for every $n\in \N$, and
$u_n\to u$ with 
 $\sup_n \calE(u_n)<\infty$. Gradient flows and generalized gradient systems featuring this subdifferential notion were analyzed in \cite{RossiSavare06} and \cite{MiRoSa13NADN}, respectively.  

 \medskip

 For the minimum problem \eqref{min-Banach}, we will stick to notation
\eqref{resolvent-set}, namely
\begin{itemize}
\item[\textbullet] $\phi= \phi(\sigma)$ for the value functional
  associated with the above minimum problem; in fact, we shall also refer to
  $\phi(\ugi;\cdot)$ as \emph{marginal function}. We remark for later use that,
  in analogy to the metric case in \cite{AmGiSa05GFMS}, it was proved in
  \cite[Lem.\,6.1]{MiRoSa13NADN} that
  \begin{equation}
    \label{added-last-mom}
    \lim_{\sigma \downarrow 0}  \phi(\sigma) = \calE(\ugi)\,,
  \end{equation} 
  where the superlinearity of $\calR$ in \eqref{Diss} is essential. 

\item[\textbullet] and $J_\sigma$ for the set of minimizers.
\end{itemize}
Additionally, in the Banach setup under Hypothesis \ref{h:X}  every
$\wt u_\sigma \in J_\sigma$ satisfies the Euler-Lagrange equation for
\eqref{min-Banach}, namely
\begin{equation}
\label{EL-BAN}
0 \in \pl\calR\big(\frac1\sigma( \wt u_\sigma{-}\ugi)\big) + \wt\xi_\sigma \quad
\text{ with } \ \wt\xi_\sigma \in \pl\calE(\wt u_\sigma).
\end{equation}
  This follows as in  \cite[Prop.\,4.2]{MiRoSa13NADN}
 (which, in turn, relies on results from \cite{Mordu-I}),  
 whose proof also uses that 
$\calR$ is locally Lipschitz around the point $\frac1\sigma( \wt
u_\sigma{-}\ugi)$;
note that this is granted by Lemma \ref{le:calR*}.
In fact,  \cite[Prop.\,4.2]{MiRoSa13NADN} would apply also if we worked with the notion of 
\emph{limiting subdifferential} previously mentioned.

 

 Furthermore, under Hypothesis \ref{h:X} the following holds:
\begin{enumerate}
\item We have $ J_\sigma \neq \emptyset$ for all $ \ugi \in
  \dom(\calE) $ and all $ \sigma>0$.

  For this, it suffices to observe that any infimizing sequence for
  $\Phi_\sigma( \cdot) $ is bounded in $\Spx$ (thanks to the facts that
  $\calE$ is bounded from below and $\calR$ has superlinear growth), and to
  resort to the weak lower semicontinuity of $\calE$ and $\calR$.
\item The multivalued mapping
  $\Jugi    \colon (0,\infty) \rightrightarrows \Spx$;
  $\sigma \mapsto \Jugi(\sigma):= J_\sigma$
 (cf.\ \eqref{Gamma-mapping})  
   is upper semicontinuous from $\R$ to
  $\Spx$, in the sense that inclusion \eqref{multivalued-usc-1} holds (even
  for the Kuratowski upper limit
  $ \Ls_{n\to\infty}^{\mathrm{weak}} \Jugi(\sigma_n)$ defined in terms of the
  \emph{weak} topology on $\Spx$).  Therefore, by
  \cite[Thm.\,3.22]{Castaing-Valadier77} we may conclude the existence of a
  (strongly)  \emph{measurable} selection
  $(0,\infty) \ni \sigma \to \wt{u}_\sigma \in J_\sigma $.
\end{enumerate} 

We conclude this section with the example of a $\GBGS$ $(\Spx, \calE, \calR)$
fulfilling the conditions from Hypothesis \ref{h:X}. To keep the exposition
simple, we confine the discussion to a G\^ateaux differentiable energy 
$\calE$, where $\pl\calE(u)$ is always a singleton. This is not a significant
restriction  when we revisit Example \ref{PDE-example} later on, 
because our focus will rather be on the properties of the dissipation
potential $\calR$. Nonetheless, it would not be difficult to adjust the
conditions in such a way as to allow for a nonsmooth, but $\lambda$-convex,
potential $W$ in \eqref{ex-diss-E} below.

\begin{example}  \upshape
\label{PDE-example}
We consider  a $\GBGS$ $(\Spx, \calE, \calR)$ such that 
\begin{itemize}
\item[\textbullet] 
$\Spx = \rmL^p(\Omega)$, $p>1$, with
  a bounded Lipschitz domain $\Omega\subset\R^d$;
\item[\textbullet] the dissipation potential $\calR: \rmL^p(\Omega) \to [0,\infty)$
 is  defined by $ \calR(v) = \int_\Omega \rmR(v(x)) \dd x $ with
  \begin{equation}
  \label{ex-diss-R}
  \begin{gathered}
  \rmR: \R \to [0,\infty) \text{ convex, s.t. } \rmR(0)=0 \text{ and }
  \\
    \exists\, \kappa, \, K >0\, \ \forall\, x, y \in \R\, : \ 
   \begin{cases}
   \rmR(x) \geq \kappa |x|^p - K,
   \\
   \rmR^*(y) \geq \kappa |y|^{p'} - K, 
  \end{cases}
  \end{gathered}
  \end{equation}
  where $p'$ is the dual exponent to $p$.
\item[\textbullet] the energy functional
  $\calE:\rmL^{p}(\Omega)\to(-\infty,\infty]$ features a (possibly) nonconvex,
  but lower order, perturbation of the Dirichlet integral, i.e.\ it is defined
  via
  \begin{subequations}
    \label{example-energy}
    \begin{equation}
      \label{ex-diss-E}
      \calE(u) : =\left\{
        \begin{array}{cl}
          \! \int_{\Omega}\frac{1}{2}|\nabla u|^{2} {+} W(u)\dd x & \text{for }
          u\in\rmH_{0}^{1}(\Omega)\text{ and }W(u)\in\rmL^{1}(\Omega)\,,
          \\
          \infty & \text{otherwise}\,.
        \end{array}
      \right.
    \end{equation}
    Along the footsteps of \cite[Sec.\,7]{RMS08}, we require for the potential
    energy density $W: \R \to \R $ that $ W\in\mathrm{C}^{2}(\R) $ and
    \begin{equation}
      \label{ass-W}
      \begin{aligned}
        \exists\, C_W>0 \ \exists\, s_p\in(1,\tfrac{p_d}{p'}) \ \forall \,
        r\in\R:\ \ \left\{
          \begin{array}{@{}ll}
            &  W''(r)\geq-C_W\,,
            \\
            &  W(r)\geq-C_W\,,
            \\
            &  |W'(r)|\leq C_W (1{+} |r|^{s_{p}})\,.
          \end{array}
        \right.
      \end{aligned}
    \end{equation}
  \end{subequations}
\end{itemize}
As shown in \cite[Sec.\,7]{RMS08} (cf.\ also 
\cite[Sec.\,4.2]{Mielke-Rossi-Stephan}), $\calE$ complies with the lower
semicontinuity and closedness conditions \eqref{HypEneBan},  where 
$\pl\calE(u)= \{\, {-} \Delta u {+} W'(u)\,\} \subset \rmH^{-1}(\Omega)=
\rmH^1_0(\Omega)^*$.
\end{example}

\subsection{The De Giorgi estimate and identity in the Banach setup}
\label{ss:Banach-statements}

In order to state the Banach-space versions of the metric energy identity
\eqref{MEIsec2} and of  the De Giorgi estimate
\eqref{eq:MetrDeGiEstimate}, we are naturally led to introduce the following
object, which corresponds to the ``slope part of the dissipation'' in
\eqref{eq:MetrDeGiEstimate}.

\begin{definition}[$\calR$-slope of the energy]
  Let $u \in \dom(\pl \calE)$. The quantity
  \begin{equation}
    \label{rslope}
    \rslo u := \inf\bigset{\calR^*({-}\xi)}{ \xi\in \pl\calE(u)} 
  \end{equation}
  is called \emph{$\calR$-slope} of the energy functional $\calE$ at $u$.
\end{definition}
 Indeed, as a consequence of the coercivity property
\eqref{coercivityRstar}  of $\calR^*$, guaranteeing that infimizing
sequences for the above minimum problem are bounded in $\Spx^*$, and of the
closedness property \eqref{HypEneBan-2}, the infimum in \eqref{rslope} is
attained, i.e.\ we have that
\begin{subequations}
  \label{properties-R-slope}
\begin{equation}
  \label{min-rslope-attained}
  \rslo u = \min\bigset{\calR^*({-}\xi)}{ \xi\in \pl\calE(u)} \,. 
\end{equation}
 To see this, consider $\xi_n \in \pl \calE(u)$ with
$\calR^*(-\xi_n)\to \rslo u$, then the coercivity of $\calR^*$ allows to
extract a subsequence (not relabeled) such that $\xi_n\weak \xi$ and the
sequence $(u,\xi_n) \weak (u,\xi)$ satisfies the assumptions of
\eqref{HypEneBan-2} with energy level $E=\calE(u)$. Hence, we have
$\xi\in \pl\calE(u)$ and $\rslo u=\calR^*(-\xi)$ follows.  Likewise, it is
immediate to check that, for every $(u_n)_n,\, u \in \Spx$  in some energy
sublevel (recall that the closedness property \eqref{HypEneBan-2} holds under
an energy bound), we have  
\begin{equation}
  \label{lsc-slope}
 (u_n)_n, \, u \in S_E, \ \       u_n \rightharpoonup u \ \
  \Longrightarrow \ \  \liminf_{n\to\infty} \rslo{u_n} \geq \rslo u\,.
\end{equation}
\end{subequations}

Obviously, one may expect that the $\calR$-slope $\rslo{ \wt{u}_\sigma}$,
evaluated along a (measurable) selection
$\sigma \mapsto \wt{u}_\sigma \in J_\sigma  $, will play the role of the
term $\psi^*\big(\big|\pl\calE \big|(\wt u_\sigma) \big)$ in the Banach-space
version of \eqref{eq:MetrDeGiEstimate}.  However, a more careful comparison
with the metric identity \eqref{MEIsec2}, featuring the quantity
\[
  \tfrac1\sigma \mathsf{d}_\sigma^-
  = \inf \bigset{ \| \tfrac{1}\sigma
    (u{-}\ugi)\|}{ u \in J_\sigma  },
\]
suggests that the notion of slope has to be adjusted. In fact, it is expedient
to bring into the picture the additional structure available in the Banach
setup, namely the fact that every $ u \in J_\sigma $ fulfills the
Euler-Lagrange equation \eqref{EL-BAN}.

We thus set forth a ``\emph{conditioned} slope part of the dissipation'',
defined by the minimization of the dual dissipation potential $\calR^*$ over
selections $\xi \in \partial\calE(u) $ that \emph{additionally} satisfy the
Euler-Lagrange equation. Accordingly, we introduce a multivalued operator
which encodes the validity of \eqref{EL-BAN}. We shall refer to these two
objects as \emph{conditioned $\calR$-slope} and \emph{conditioned
  subdifferential}, respectively.  Notice that the conditioned
$\calR$-slope is naturally defined on the graph of the multivalued operator
$\Jugi$, namely on the set
\[
\mathrm{Graph} (\Jugi) = \{ (\sigma,u) \in  (0,\infty){\times}X \,: 
\ u \in \Jugi(\sigma) = J_\sigma \}\,.
\] 

\begin{definition}[\emph{Conditioned} subdifferential\,/\,slope of energy]
\label{def:CondSubdSlope} \hfill The multivalued mapping 
$ \ELname \calE: \mathrm{Graph} (\Jugi) \rightrightarrows \Spx^*$  defined
by
\begin{equation}
  \label{EL-slope}
  \EL {\sigma}{u}:=\Bigset{\xi \in  \Spx^*}{\xi\in \pl\calE(u) 
    \text{ and } {-}\xi \in \partial\calR\big( \frac1\sigma(u{-}\ugi)\big)}
\end{equation}
is called \emph{\nameop of the energy}  $\calE$. The quantity
\begin{equation}
  \label{eq:ConditSlope}
  \wrslo u{\sigma} :=\inf\bigset{\calR^*({-}\xi)}{ \xi\in  \EL {\sigma}{u}}  \quad \text{for all } (\sigma,u) \in \mathrm{Graph} (\Jugi)  
\end{equation} 
 is called \emph{\nameslope} of $\calE$ at $u \in \dom(\pl \calE)$.
\end{definition}

Although with slight abuse, the notation for $ \ELname \calE$ highlights the
geometry induced by the dissipation potential through the Euler-Lagrange
equation \eqref{EL-BAN}.

Clearly, since every $u \in J_\sigma$ fulfills the Euler-Lagrange
equation \eqref{EL-BAN}, we have $ \EL {\sigma}{u} \neq \emptyset$ and
thus $ \wrslo u{\sigma} < \infty$  for all $(\sigma,u) \in \mathrm{Graph} (\Jugi)$.  Obviously, in general there holds
\begin{equation}
\label{FF-later}
  \wrslo u{\sigma} \geq \rslo u \qquad \text{for all }( \sigma, u) 
  \in \mathrm{Graph} (\Jugi),
\end{equation}
which seems to suggest $ \wrsloname $ as the `right' object for the attainment
of an equality in  the  De Giorgi estimate, cf.\ also Example \ref{ex:2-revisited}
below.

We are now in the position to precisely introduce the estimate/identity whose
validity we are are going to address hereafter.

\begin{definition}
\label{def:MEI-DG-Banach}
We say that a (strongly) measurable selection
$(0,\infty) \ni \sigma \mapsto \wt{u}_\sigma \in J_\sigma  $ fulfills, on
some interval $(0,\sigma_*)$,
\begin{itemize}
 \item[the \textbf{simple De Giorgi estimate}]  if for all
  $\sigma \in (0,\sigma_*)$ we have 
  \begin{equation}
    \label{eq:BGS.DGE.Simple} {
      \calE(\wt u_\sigma) + \sigma \calR\big(\frac1\sigma (\wt u_\sigma {-}
      \ugi)\big) + \int_{\rho=0}^\sigma \rslo { \wt{u}_\rho } \dd \rho
      \leq \calE(\ugi);   }
  \end{equation}  
 \item[the \textbf{improved De Giorgi estimate}] if for all
  $\sigma \in (0,\sigma_*)$ we have 
  \begin{equation}
    \label{eq:BGS.DeGiEstim}
    \begin{aligned}
      \calE(\wt u_\sigma) + \sigma \calR\big(\frac1\sigma (\wt u_\sigma {-}
      \ugi)\big) + \int_{\rho=0}^\sigma \wrslo { \wt{u}_\rho } \rho \dd \rho
      \leq \calE(\ugi);
    \end{aligned}
  \end{equation} 
\item[the \textbf{De Giorgi identity}] if \eqref{eq:BGS.DeGiEstim} holds as an
  equality, i.e.\ for all $\sigma \in (0,\sigma_*)$:
  \begin{equation}
    \label{MEI-Ban-cond}
    \calE(\wt u_\sigma) + \sigma \calR\big(\frac1\sigma (\wt
    u_\sigma {-} \ugi)\big)  + \int_{\rho=0}^\sigma     \wrslo  {\wt{u}_\rho}
    {\rho}  \dd  \rho =  \calE(\ugi)\,.
  \end{equation}
\end{itemize}
\end{definition}

Later on, in Section \ref{ss:4.2} we will provide a sufficient condition for
the attainment of the infimum in the definition
\eqref{eq:ConditSlope} of $\wrsloname$, which is also related to the
validity of the analog of the lower semicontinuity properties
\eqref{lsc-slope}.  This will lead to the existence of measurable selections
\[
  (0,\infty) \ni \sigma \mapsto \wt{\xi}_\sigma \in \EL
  {\sigma}{\wt{u}_\sigma}, \quad \text{with $\sigma \mapsto\wt{u}_\sigma \in \Jugi (\sigma)$
    a measurable selection,}
\]
such that estimate \eqref{eq:BGS.DeGiEstim} rephrases as
\[
    \calE(\wt u_\sigma)  + \sigma \calR\big(\frac1\sigma 
    (\wt u_\sigma {-} \ugi)\big)
 + \int_0^\sigma   \calR^*({-}\wt \xi_\rho)
 \dd \rho  \leq  \calE(\ugi)\,,
\]
and analogously for \eqref{MEI-Ban-cond}. 

We conclude this section by using Example \ref{ex:SlopeNotContin}, revisited in
the Banach setup, to convey the idea  that the finer information
encoded in the {\nameslope} of the energy may play a key role for the
attainment of  the  De Giorgi identity.
 
\begin{example}[Example \ref{ex:SlopeNotContin} in the Banach setup]
\label{ex:2-revisited} \slshape
We now treat Example \ref{ex:SlopeNotContin}  as
a $\GBGS$  in the Hilbert or Banach space $X=\R$: 
\[
M=\R,\quad \calE(u)=\max\{u,0\} =: u^+ ,\quad \calR(v) = \frac12 v^2.
\]
Since $\calE$ is convex, its Fr\'echet subdifferential
$\pl \calE: \R \rightrightarrows \R$ coincides with the subdifferential  in the sense of convex
analysis, and it is set-valued with $\pl\calE(0)=[0,1]$. The $\calR$-slope of
the energy is given by
\[
\rslo u=\frac12 \text{ for }u>0 \quad \text{ and } \quad 
 \rslo u=0 \text{ for }u\leq 0 .
\]
\par
As before, starting from
$\ugi=1$
 we calculate the variational interpolant  $\wt u_\sigma = (1{-}\sigma)^+=\max\{ 1{-}\sigma, 0\}$, so that 
\[
\pl \calE(\wt u_\sigma) = \begin{cases}
\{1\} &\text{ if } \sigma \in (0,1),
\\
[0,1] & \text{ if } \sigma \geq 1.
\end{cases}
\] 
 Therefore,  for $\sigma >0$
 we have 
 \[
   \calE(\ugi)-   \calE(\wt u_\sigma)  - \sigma \calR\big(\tfrac1\sigma (\wt
u_\sigma {-} \ugi)\big) = 1 {-}  (1{-}\sigma)^+- \tfrac1{2\sigma} |(1{-}\sigma)^+{-}1|^2
 = \begin{cases}
\frac1{2}\sigma &\text{if } \sigma \leq 1,
\\
1{-} \frac1{2\sigma} & \text{if } \sigma \geq 1.
\end{cases}
\]
Hence, the  estimate
\begin{equation}
  \label{concrete-DGE}
  \int_0^\sigma \calR^*\big({-}\wt\xi_\rho\big) \dd \rho \leq \calE(\ugi)
   -   \calE(\wt u_\sigma)  - \sigma \calR\big(\frac1\sigma (\wt
  u_\sigma {-} \ugi)\big) 
\quad \text{with } \wt\xi_\rho   \in \pl \calE(\wt u_\rho)
\end{equation}
holds for $\sigma \in (0,1]$, while for $\sigma > 1$ it only holds for
selections $\sigma \mapsto \wt\xi_\sigma \in \pl \calE (\wt u_\sigma) $
with
\[
  \int_0^\sigma \calR^*\big({-}\wt\xi_\rho \big) \dd \rho = \int_0^1 \frac12 \dd
  \rho + \int_1^{\sigma} \frac12 \big|{-}\wt\xi_\rho \big|^2 \dd \rho \leq 1 -
  \frac{1}{2\sigma}\quad \Leftrightarrow \quad 
  \int_1^\sigma \big|\wt\xi_\rho\big|^2 \dd \rho \leq 1-\frac1\sigma. 
\]
Since $ \rslo {\wt u_\rho}=0$ for $\rho \geq 1$, \eqref{concrete-DGE} is
certainly satisfied if we replace the integrand $\calR^*\big({-}\wt\xi_s\big)$
by $ \rslo {\wt u_s}$;  thus the simple De  Giorgi estimate
\eqref{eq:BGS.DGE.Simple} is clearly satisfied. 

However,   recalling that for $\sigma \geq 1$ we have
$\pl \calE(\wt u_\sigma) =[0,1]$, a wrong choice of
$\wt\xi_\sigma \in \pl \calE(\wt u_\sigma) $ can violate \eqref{concrete-DGE}.

Finally,  we consider the Euler-Lagrange equation \eqref{EL-BAN} in the
current Banach setting:
\[
0 \in \wt\xi_\sigma + \pl\calR\big(\frac1\sigma(\wt
u_\sigma{-}\ugi)\big) \quad \text{with } \wt\xi_\sigma \in \pl\calE(\wt u_\sigma) .
\]
Since $\calR$ is smooth with $\pl\calR(v)=\{v\}$, there is a unique solution,
namely 
\[
  \wt\xi_\sigma = -\frac1{\sigma} (\wt u_\sigma - \ugi) = \begin{cases} 1 &
    \text{for }\sigma \in (0,1]\,,\\ 1/\sigma & \text{for } \sigma \geq  1\,.
  \end{cases}
\]
In particular,  we have 
\begin{align*}
\text{for  } \sigma \in (0,1]:&\quad  \EL {\sigma}{\wt u_\sigma}=\{1\} 
\quad \text{and} \quad \calS_\calR(\wt
  u_\sigma)= \calC_\calR(\sigma; \wt u_\sigma) = \frac12,
\\
\text{for } \sigma >1:&\quad   \EL {\sigma}{\wt
  u_\sigma}=\big\{\frac1\sigma \big\} \quad \text{and} 
  \quad 0=\calS_\calR(\wt u_\sigma) \lneqq \calC_\calR(\sigma; \wt
  u_\sigma)=\frac1{2\sigma^2}\:.
\end{align*}
Moreover,  equality in \eqref{concrete-DGE} can be achieved only 
by choosing the special selection 
\[
\wt\xi_\sigma \in    \EL {\sigma}{\wt u_\sigma} \qquad \foraa\ \sigma >0\,.
\]
Thus, the `good selections' of $\pl \calE (\wt u_\sigma)$ for the purpose of De
Giorgi identity is prescribed by the Euler-Lagrange equation
\eqref{EL-BAN}.  

In summary, we have the simple De Giorgi estimate \eqref{eq:BGS.DGE.Simple} for
$\rsloname$ (with strict inequality for $\sigma>1$), and we have  the 
De Giorgi identity \eqref{MEI-Ban-cond} for $\wrsloname$.
\end{example}

\subsection{Radially differentiable potentials}
\label{ss:radial-diff}

This section revolves around a structural property for dissipation potentials
that will be crucial for obtaining   the  De Giorgi identity. 

\begin{definition}
\label{def:radial-differentiability}
We say that a dissipation potential $\calR: X 
\to [0,\infty] $ is \emph{radially differentiable} if 
\begin{equation}
  \label{eq:BGS.StraCond2}
  \forall\, v\in   \dom(\calR)  : \quad  (0,1) 
  \ni \lambda \mapsto f(v;\lambda):=\calR(\lambda v)\text{ is
    differentiable. } 
\end{equation}
\end{definition}
 Observe that, if $\mathrm{dom}(\calR)$ is \emph{open}, then  \eqref{eq:BGS.StraCond2} in fact implies that 
$\lambda \mapsto f(v;\lambda)$ is well-defined and differentiable on the interval $(0,1+\epsilon)$ for some $\epsilon>0$.  
\par
Clearly, both differentiability on the one hand, and positive homogeneity
(i.e., $\calR(\lambda v) = \lambda^p \calR(v)$ for some $p \geq 1$) on the
other, are sufficient conditions for \eqref{eq:BGS.StraCond2}.  Thus, linear
combinations of convex differentiable, or positively homogeneous,
potentials are radially differentiable.

We are now going to show that radial differentiability is equivalent to another
structural property that was used in \cite{MiRoSa13NADN} to prove  the 
De Giorgi estimate \eqref{eq:BGS.DeGiEstim}.

\begin{proposition}
\label{prop:equivalence-structural}
 A dissipation potential $\calR$ is radially differentiable if and only if
 \begin{equation}
  \label{eq:BGS.StrangeCond}
  \xi_1,\xi_2 \in \pl \calR(v) \quad \Longrightarrow \quad
  \calR^*(\xi_1)=\calR^*(\xi_2). 
\end{equation}
 \end{proposition}
 In view of the well-known convex-analysis relation 
\[
 \calR(v)  +  \calR^*(\xi) =  \langle \xi,v\rangle  \qquad 
 \text{for all } \xi \in \partial  \calR(v),
\]
condition \eqref{eq:BGS.StrangeCond} is, in turn, equivalent to the property that 
\[
  \xi_1,\xi_2 \in \pl \calR(v) \quad \Longrightarrow \quad \langle \xi_1,v\rangle = \langle \xi_2,v\rangle.
\]
This ensures that the (a priori multivalued) mapping   (whose domain is clearly given by $ \dom(\partial\calR)$)  
\begin{equation}
\label{equal-result}
\Prod: \Spx \rightrightarrows [0,\infty),\qquad  \Prod(v) := \langle \xi,v\rangle   \text{ for all } 
\xi \in \partial  \calR(v), \quad \text{ is single-valued}
\end{equation}
which will prove useful in Section \ref{ss:regularity-added}
ahead. 
\par 
Let us now address the proof of 
 the equivalence of \eqref{eq:BGS.StrangeCond} and \eqref{eq:BGS.StraCond2}: We start by observing that 
for every  $v\in  \dom(\calR)  $   the mapping $f(v;\cdot)$ is convex.   Its subdifferential
in the sense of convex analysis  $\partial f(v;\cdot): (0,\infty) \rightrightarrows \R$ is given by 
\begin{subequations}
\label{subdiff-f}
\begin{equation}
\label{subdiff-f-1}
\partial f(v;\lambda) = [ f'_-(v;\lambda), f'_+(v;\lambda)] \qquad \text{for all } \lambda >0, \ v \in   \dom(\calR) \,, 
\end{equation}
where $ f'_\pm(v;\cdot)$ are 
 the one-sided derivatives of  the mapping
$ f(v;\cdot)$, 
\begin{equation}
\label{subdiff-f-2}
\begin{aligned}
  & f'_+(v;\lambda):=\lim_{h\to 0^+} \frac1h\big( f(v;\lambda{+}h)-
  f(v;\lambda)\big),
  \\
  & f'_-(v;\lambda):=\lim_{h\to 0^+} \frac1h \big( f(v;\lambda)-
  f(v;\lambda{-}h)\big) \,,
\end{aligned}
\end{equation}
(the above limits exist by convexity of $f(v;\cdot)$).
\end{subequations}
Then, we have the following result.

\begin{lemma}
\label{l:radial-derivatives}
For every   $v \in 
 \dom(\calR)$   and $\lambda>0$
\begin{subequations}
\label{eq:LeftRightDeri}
\begin{align}
\label{eq:LeftRightDeri+}
  f'_+(v;\lambda) 
  & =
  \max_{\xi\in \pl\calR(\lambda v)} \langle \xi,v\rangle = \calR(\lambda v) - \min_{\xi\in
    \pl\calR(\lambda v)}\calR^*(\xi), 
\\ 
\label{eq:LeftRightDeri-}
  f'_-(v;\lambda) 
  & =
  \min_{\xi\in \pl\calR(\lambda v)} \langle \xi,v\rangle = \calR(\lambda v) - \max_{\xi\in
    \pl\calR(\lambda v)}\calR^*(\xi).
\end{align}
\end{subequations}
\end{lemma}
\begin{proof}
It can be easily checked that, for any   $v\in  \dom(\calR)$   and $\lambda>0$, there holds
\[
\ell  \in \partial f(v;\lambda) \   \text{ if and only if }  \ \exists\, \xi \in \partial \calR(\lambda v) \text{ s.t. } \ell =\pairing{}{}{\xi}{v}\,.
\]
Combining this with \eqref{subdiff-f-1}, we immediately deduce
\eqref{eq:LeftRightDeri}.  
\end{proof}

We are now in a position to carry out the
\\[0.4em]
\begin{proof}[Proof of Proposition \ref{prop:equivalence-structural}] The mapping
$f(v;\cdot)$ from \eqref{eq:BGS.StraCond2} is differentiable at $\lambda =1$ if
and only if $ f'_-(v;1) =f'_+(v;1) $. This is in turn equivalent, by
\eqref{eq:LeftRightDeri}, to the fact that
$ \min_{\xi\in \pl\calR(v)}\calR^*(\xi) = \max_{\xi\in
  \pl\calR(v)}\calR^*(\xi)$, i.e.\ \eqref{eq:BGS.StraCond2}.  
\end{proof}

\begin{remark}\slshape
\label{rmk:more-on-radially}
Let us focus on `metric-like' dissipation potentials of the form
\begin{equation}
\label{metric/like}
\calR(v):= \psi(\| v \|_X) \qquad \text{for every } v \in X 
\end{equation}
with $\psi: [0,\infty) \to [0,\infty) $ convex with superlinear growth at
infinity, as in Definition \ref{def:MGS}.  For the associated mapping
$f(v;\cdot)$, $v \in \Spx$ (note that $ \dom(\calR)=X$, in this case), 
 we have
\[
\partial f(v;\lambda) = \|v\|_X \partial \psi (\| \lambda v \|_X) \qquad
\text{for all }\lambda>0\,. 
\]
Therefore, $\calR$ is radially differentiable in the sense of
\eqref{eq:BGS.StraCond2} if and only if $\partial \psi$ is single-valued. We
thus retrieve the smoothness requirement on $\psi$ from Hypothesis \ref{h:M}.
\end{remark}

\begin{remark}[Example \ref{PDE-example} revisited. ]
\label{rmk:ex-revisited}
{\slshape Let us get back to the dissipation potential
$\calR:\rmL^p(\Omega) \to [0,\infty)$, $p\geq 1$, from Example
\ref{PDE-example}, i.e.\ $\calR(v) = \int_\Omega \mathrm{R}(v(x)) \dd x $, with
the dissipation density $\mathrm{R} : \R \to [0,\infty)$ satisfying conditions
\eqref{ex-diss-R} (so that $ \dom(\calR)= \rmL^p(\Omega) $).    In that setting, for a given
$v\in \dom(\partial\calR)$ we have that
\[
\xi \in \partial \calR(v)  \quad \Longleftrightarrow \quad \xi(x) \in
\partial\mathrm{R}(v(x)) \ \  \foraa\, x \in \Omega\,. 
\]
Therefore, it is natural to address the relations between radial
differentiability of $\calR$ and radial differentiability of $\mathrm{R}$,
which is clearly equivalent to differentiability of $\mathrm{R}$ in
$\R{\setminus}\{0\}$.  We now check that
\begin{equation}
\label{equivalence-radial-differentiabilities}
\calR \text{ is radially differentiable } \ \Longleftrightarrow   \ \mathrm{R}
\text{ is differentiable in }  \R{\setminus}\{0\}. 
\end{equation}

Indeed, suppose that $\mathrm{R}: \R \to [0,\infty)$ is radially
differentiable.
Then, in view of the characterization provided by Prop.\
\ref{prop:equivalence-structural}, for all $v\in \dom(\partial\calR)$
and $\xi_1,\, \xi_2 \in \partial \calR(v)$ we have
\[
  \begin{aligned}
    & \mathrm{R}^*(\xi_1(x)) = \mathrm{R}^*(\xi_2(x)) \ \ \foraa\, x \in \Omega
    \quad \text{ and thus }
    \\
    & \calR^*(\xi_1) = \int_\Omega \mathrm{R}^*(\xi_1(x)) \dd x =\int_\Omega
    \mathrm{R}^*(\xi_2(x)) \dd x = \calR^*(\xi_2) \,,
  \end{aligned}
\]
hence $\calR$ is radially differentiable.

The converse implication holds since, if $\mathrm{R}$ is not radially
differentiable, then there exist $\hat v, \hat{\xi}_1, \hat{\xi}_2 \in \R$ such
that $\hat{\xi}_i \in \partial\mathrm{R}(\hat{v})$, $i=1,2$, and
$\mathrm{R}^*(\hat\xi_1)\neq \mathrm{R}^*(\hat\xi_2)$. Then, defining
$v(x) \equiv \hat v$ and $\xi_i(x) \equiv \hat\xi_i$ for almost all
$x\in \Omega$, we obtain a triple
$(v,\xi_1,\xi_2) \in \rmL^p(\Omega)\ti \rmL^{p'}(\Omega)\ti \rmL^{p'}(\Omega)$
for which \eqref{eq:BGS.StrangeCond} fails to hold.

For later use, we point out that, if $\calR$ is radially differentiable, then
the single-valued mapping $\Prod: \rmL^p(\Omega) \to [0,\infty)$ from
\eqref{equal-result} is indeed given by
\begin{equation}
  \label{Prod-example}
  \Prod(v) = \int_\Omega \mathrm{P}(v(x)) \dd x \qquad \text{with } 
  \mathrm{P}(r) : =  \begin{cases}
    r \, \mathrm{R}'(r) &\quad\text{if }r \in \R{\setminus}\{0\},
    \\  0 &\quad\text{if } r=0\,.
  \end{cases}
\end{equation}
} 
\end{remark}

Ultimately,  as a straightforward consequence of Prop.\
\ref{prop:equivalence-structural} \ we have the following result.

\begin{corollary}
Suppose that 
 $\calR: X 
\to [0,\infty) $ is \emph{radially differentiable}. Then, 
we have that 
  \[
   \wrslo u{\sigma}=  \calR^*({-}\xi)   \qquad \text{for all } \, 
   \xi \in \EL {\sigma}{\wt{u}_\sigma} \text{ and } u 
 \in  J_\sigma\,.
  \] 
\end{corollary}

In particular, as soon as a (measurable) selection $\sigma \mapsto \wt u_\sigma
\in J_\sigma $ fulfills the  De Giorgi estimate/identity, then
\eqref{eq:BGS.DeGiEstim}/\eqref{MEI-Ban-cond} hold with the slope part of the
time integrated dissipation,  namely $ \int_0^\sigma \wrslo { \wt{u}_\rho
}{\rho} \dd \rho$, given by $ \int_0^\sigma \calR^*({-}\wt\xi_\rho) \dd
\rho$ for \emph{any} (measurable) selection $\sigma \mapsto \wt\xi_\sigma \in
\EL {\sigma}{\wt{u}_\sigma}$.  Furthermore, if the  De Giorgi identity
holds along an interval $(0,\sigma)$, the corresponding \emph{optimal} 
integrand $ \wrslo { \wt{u}_\rho }{\rho} $ fulfills
 \begin{equation}
 \label{optimal-value}
  \wrslo { \wt{u}_\rho }{\rho} 
 =  \frac{\rmd}{\rmd s} \big( s
\calR\big(\frac1s(\wt u_\rho{-}\ugi)\big)\big)|_{s=\rho } 
\qquad \foraa\, \rho \in (0,\sigma)\,.
\end{equation}

In Theorem  \ref{th:DGL.RadialDiff} ahead we will indeed prove  the 
De Giorgi identity under the condition that the dissipation potential is
radially differentiable. 

\subsection{Tools}
\label{ss:3.4}

In this section we collect some preliminary results that will be used in the
proof of Theorem~\ref{th:DGL.RadialDiff}.

\subsubsection*{Radial derivatives of dissipation potentials} 
Our first result is in the same spirit of Lemma \ref{l:radial-derivatives}, in
that it provides some information on the left and right derivatives of the
function $g(t)=t \calR\big(\frac1t\,v\big)$, $t>0$, which also features in
\eqref{optimal-value}.

\begin{lemma}[Radial derivative]
\label{le:RadiDeri}
For fixed   $v\in \mathrm{dom}(\calR)$   define the function $g(t)=t
\calR\big(\frac1t\,v\big)$. Then, $g$ is convex, decreasing, and satisfies
$g \in \rmC^\mafo{Lip}_\mafo{loc}({]0,\infty[})$. For all $t>0$ the left
derivative $g'_-(t)$ and the right derivative $g'_+(t)$ exist, are
non-decreasing, and satisfy
\begin{align*}
g'_-(t)&:=\lim_{h\to 0+}\frac1h\big(g(t)-g(t{-}h)\big) =
- \max\bigset{ \calR^*(\eta)}{ \eta\in \pl\calR(\frac1t v)}
\\
&\leq - \min\bigset{ \calR^*(\eta)}{ \eta\in \pl\calR(\frac1t v)} =
\lim_{h\to 0+}\frac1h\big(g(t{+}h)-g(t)\big) =: g'_+(t) <0. 
\end{align*}
Moreover, $t\mapsto g'_-(t)$ is continuous from the left and $t\mapsto
g'_+(t)$ is continuous from the right. 
\end{lemma}

In particular, if $\calR$ is radially differentiable, then $g$ is
continuously differentiable with $g'(t)=g'_\pm(t)=- \calR^*(\xi)$ for all $\xi
\in \pl\calR(\frac1t\,v)$, in accordance with Lemma \ref{l:radial-derivatives}. 

\subsubsection*{Fundamental lemma for marginal functions}
 Recall that $(0,\infty) \ni \sigma \mapsto \phi(\sigma)$  is   the marginal function  
of $(\sigma,u)\mapsto \Phi_\sigma(
u)$.  Our next result shows that 
 one-sided differentiability of $\Phi$ with respect to
$\sigma$ provides bounds on the one-sided derivatives of $\phi$, for which  the behavior of
$\Phi$ with respect to $u$ is not really important. In the following result, in fact,  we do not use the
special form of $\Phi$, but only its left and right differentiability with
respect to $\sigma$. 

\begin{proposition}[Derivatives of marginal functions]
\label{prop:DeriMarg} 
We have $\phi\in \rmC^\mafo{Lip}_\mafo{loc}({]0,\infty[})$ and for all $\sigma
>0$ the following estimates hold
\begin{subequations}
 \label{eq:BGS.MargiEst}
\begin{align}
  \label{eq:BGS.MargiEst.a}
&\liminf_{h\to 0^+} \frac1h \big(\phi(\sigma)- \phi(\sigma{-}h)\big)  \geq 
\sup\bigset{\rmD^-_\sigma \Phi_\sigma(w)}{ w \in J_\sigma}=: 
 \delta^-_\sigma, \\
 \label{eq:BGS.MargiEst.b}
& \limsup_{h\to 0^+} \frac1h \big(\phi(\sigma{+}h)- \phi(\sigma)\big) \leq
 \inf\bigset{\rmD^+_\sigma \Phi_\sigma(w)}{ w \in J_\sigma}=:
\delta^+_\sigma
\end{align}
\end{subequations}
  with $\rmD^\pm_\sigma \Phi_\sigma(w)$ the one-sided derivatives of the map $\sigma \mapsto \Phi_\sigma(w)$. 
\end{proposition}
\begin{proof}
Consider $0 < r < \sigma$, then by the marginal property of $\phi$ we have
\begin{align*}
\Phi_\sigma(\wt u_r) - \Phi_r(\wt u_r) \geq 
\Phi_\sigma(\wt u_\sigma) - \Phi_r( \wt u_r)
&= \phi(\sigma) - \phi(r)
\geq \Phi_\sigma(\wt u_\sigma) - \Phi_r(\wt u_\sigma). 
\end{align*}
Hence, the local Lipschitz property of $\phi$ follows from that of $\Phi(
\cdot,w)$. Setting $r=t{-}h$ yields
\[
\frac1h\big( \phi(\sigma)- \phi(\sigma{-}h)\big) \geq \frac1h\big(\Phi_\sigma(
\wt u_\sigma) - \Phi_{\sigma{-}h}(\wt u_\sigma) \big). 
\] 
Taking the liminf for $h\to 0^+$ first and then the supremum over $\wt u_\sigma
\in J_\sigma$ gives  \eqref{eq:BGS.MargiEst.a}. 

Similarly, we can replace $(r,\sigma)$ by $(\sigma,\sigma{+}h)$ to obtain
\[
\frac1h\big( \phi(\sigma{+}h)- \phi(\sigma)\big) \leq \frac1h\big(\Phi_{\sigma {+}h}(
\wt u_\sigma) - \Phi_\sigma(   \wt u_\sigma) \big). 
\] 
Taking the limsup $h\to 0^+$ first and then the infimum over $\wt
u_\sigma\in J_\sigma$ yields \eqref{eq:BGS.MargiEst.b}. 
\end{proof}

\section{The Banach case: results}
\label{s:4}

In the upcoming Section \ref{ss:4.1} we will show that the De Giorgi estimate
previously proved in \cite[Lemma 6.1]{MiRoSa13NADN} in the radially
differentiable case, in fact improves to an equality.   Section
\ref{ss:regularity-added} will address the regularity of the variational
interpolant.  Then, in Sec.\ \ref{ss:4.2} we will drop the radial
differentiabiilty condition and extend the validity of  the  De Giorgi
estimate to general dissipation potentials.
\subsection{Equality in  the  De Giorgi estimate for radially differentiable potentials}
\label{ss:4.1}
Now we return to the variational integrand where $\Phi$ is given in the form 
\[
\Phi_\sigma( u) = \sigma\,\calR\big( \frac1\sigma(u{-}\ugi)\big) + \calE(u).
\]
In particular, the left and right derivatives of $\sigma \mapsto \Phi_\sigma(
u ) $ exist, see Lemma \ref{le:RadiDeri}, and are independent of the energy
$\calE$. Moreover, the one-sided derivatives are ordered such that
$\rmD^-_\sigma\Phi(\sigma,w) \leq \rmD^+_\sigma \Phi(\sigma,w) < 0$ for
every $w \in J_\sigma$. 

However, because of the supremum over $\rmD^-_\sigma \Phi$ and the infimum over
$\rmD^+_\sigma \Phi$, in general one cannot expect to generate one chain of
inequalities from the two estimates in \eqref{eq:BGS.MargiEst}. Even if
$J_\sigma$ is single-valued, one still arrived at  the wrong estimate.

The only case where the marginal estimates \eqref{eq:BGS.MargiEst} are useful
 to prove the  De Giorgi identity,  is exactly when $\calR$ is
radially differentiable, cf.\ \eqref{eq:BGS.StraCond2}. Then
$\rmD^-_\sigma \Phi= \rmD^+_\sigma \Phi$ implies
$\delta^-_\sigma\geq \delta^+_\sigma$ and \eqref{eq:BGS.MargiEst} leads to
\[
\liminf_{h\to 0^+}\frac1h\big(\phi(\sigma)- \phi( \sigma{-}h) \big) \geq 
\delta^-_\sigma\geq \delta^+_\sigma \geq \limsup_{h\to 0^+}\frac1h
\big( \phi( \sigma{+}h) - \phi(\sigma)\big) .
\]
From this,  the  De Giorgi identity follows easily.

 We will also address the validity of the identity involving a measurable
selection $\wt\xi_\sigma\in \EL {\sigma}{\wt u_\sigma} $ (we will often refer
to $\wt \xi$ as a \emph{force selection}).  For the existence of such a
selection (see Lemma \ref{l:conseq-of-closedness}), it is useful to impose the
following closedness condition for the conditioned subdifferential
$ \ELname \calE : (0,\infty)\ti \Spx \rightrightarrows \Spx^*$.

\begin{hypothesis}[Closedness of $\ELname \calE$]
\label{hyp:EL-closedness}
The \nameop
$\ELname \calE: (0,\infty)\ti \Spx$ $ \rightrightarrows \Spx^*$ is
closed on energy sublevels, i.e.
\begin{align}
\label{HypEneBan-EL}
\hspace{-0,5cm}\forall\, E>0: \ \left\{ 
\begin{array}{c}
  (\sigma_n, u_n,\xi_n) \rightharpoonup (\sigma, u,\xi) \   \text{ in }  
  (0,\infty) \ti \Spx \ti \Spx^*,  \\
  u_n \in S_E  \text{ and } \xi_n \in \EL{\sigma_n}{u_n}  
 \text{ for all } n \in \N
\end{array}
\right\} \  \Longrightarrow  \ \xi \in \EL{\sigma}{u}\,.
\end{align}
\end{hypothesis}

\begin{theorem}[The De Giorgi identity with radial differentiability]
\label{th:DGL.RadialDiff} 
Consider the $\GBGS$ $(X,\calE,\calR)$ satisfying Hypothesis \ref{h:X}.
Fix $\ugi\in \mafo{dom}(\calE)$ and assume that $\calR$ 
is radially differentiable, i.e.\ \eqref{eq:BGS.StrangeCond} or
equivalently \eqref{eq:BGS.StraCond2} holds. 
Then, every measurable variational interpolant
$(0,\infty) \ni \sigma \mapsto \wt u_\sigma\in  J_\sigma$ 
fulfills  the  \emph{De Giorgi identity}
\begin{subequations}
  \label{eq:DegIdent4}
 \begin{align}
 \label{eq:DegIdent4-nosel}
  \calE(\wt u_\sigma) + \sigma \calR\big(\frac1\sigma (\wt
u_\sigma {-} \ugi)\big)
 + \int_{0}^\sigma     \wrslo  {\wt{u}_\rho} {\rho}  \dd
 \rho =  \calE(\ugi)\,. 
\end{align}
 If additionally Hypothesis \ref{hyp:EL-closedness} holds, then there 
exists a measurable force selection $\sigma \mapsto \wt\xi_\sigma\in \EL 
{\sigma}{\wt u_\sigma} $ such that the 
\emph{improved De Giorgi identity with force selection} holds:
\begin{align}
  \label{eq:BGS.DeGiIdent-sel}
  \forall \, \sigma>0: \quad \calE(\wt u_\sigma) + \sigma 
  \calR\big( \frac1\sigma(\wt u_\sigma{-}\ugi)\big)
   + \int_0^\sigma \calR^*(-\wt\xi_\rho) \dd \rho = \calE(\ugi) .
\end{align}
\end{subequations}
\end{theorem}
\begin{proof} We fix a small $r>0$ and consider the marginal function
  $[r,\sigma] \ni \rho\mapsto  \phi(\rho)$. According to Proposition 
\ref{prop:DeriMarg} $\phi$ is Lipschitz and hence differentiable almost
everywhere in $[r,\sigma]$. This implies that $\phi'_-$ and $\phi'_+$ exist
and coincide almost everywhere. 

We now fix any measurable variational interpolant $\sigma\mapsto \wt u_\sigma$.
Proposition \ref{prop:DeriMarg} and Lemma \ref{le:RadiDeri} imply, via the
radial differentiability of $ \calR$, that for almost all
$\sigma \in (0,\infty)$ there holds 
\[
\phi'(\sigma) = - \calR^*(\eta) \qquad \text{for all } \eta \in
\pl\calR\big(\frac1\sigma(\wt u_\sigma{-}\ugi)\big) \, .
\]
 Thus, 
\[
\phi'(\sigma)   =  - \wrslo  {\wt{u}_\sigma} {\sigma} 
   \quad \foraa\, \sigma \in (0,\infty) \, . 
\]

Hence, $\phi(\sigma)  = \phi(r) + \int_r^\sigma \phi'(\rho) \dd \rho$ 
can be rewritten as
\[
  \phi(\sigma)= \calE(\wt u_\sigma) + \sigma \calR\big( \frac1\sigma(\wt
  u_\sigma{-}\ugi)\big)  = \phi(\ugi;r) - \int_r^\sigma \wrslo
  {\wt{u}_\rho} {\rho} \dd \rho \, . 
\] 
The superlinearity of $\calR$ implies $\phi(r)\to \calE(\ugi)$ as
$r\to 0^+$,  see \eqref{added-last-mom}. Therefore, taking the limit
$r\to 0^+$ gives the desired De Giorgi  identity \eqref{eq:DegIdent4-nosel}. 

 The last part of the assertion immediately follows from the upcoming
Lemma \ref{l:conseq-of-closedness}. 
\end{proof}

 The following result collects two straightforward consequences of Hypothesis  \ref{hyp:EL-closedness}.
\begin{lemma}
\label{l:conseq-of-closedness}
Let the  $\GBGS$    $(X,\calE,\calR)$
satisfy  Hypotheses \ref{h:X} and \ref{hyp:EL-closedness}.  
Then,
\begin{enumerate}
\item The \emph{infimum} in the definition of $\ELname \calE$ is attained, i.e.\
\begin{equation}
\label{attainment-wrslo}
\wrslo u{\sigma} :=\min \bigset{\calR^*({-}\xi)}{ \xi\in  \EL {\sigma}{u}} 
\quad \text{for all } u \in J_\sigma;
\end{equation}
\item For every measurable selection
  $(0,\infty) \ni \sigma \mapsto \wt u_\sigma \in J_\sigma $ lying in some
  energy sublevel $ S_E$, there exists a measurable selection
  \begin{subequations}
    \label{meas-sele-ARGMIN}
    \begin{equation}
      \label{def-frakA}
      \sigma \mapsto \wt \xi_\sigma \in  \mathfrak{A}_{\calR} (\sigma, 
      \wt{u}_\sigma):= \mafo{argmin} \bigset{\calR^*({-}\xi)}{\xi\in 
        \EL\sigma{\wt{u}_\sigma}}
    \end{equation}
    such that
    \begin{equation}
      \label{measurable-selection}
      \wrslo{ {  \wt u_\sigma } }{\sigma} = \calR^*({-}\wt\xi_\sigma) 
       \quad \text{for all } \sigma >0\,.
    \end{equation}
  \end{subequations}
\end{enumerate}
\end{lemma} 
\begin{proof}
Property \eqref{attainment-wrslo} is easily checked via the \emph{direct method}
by relying on \eqref{HypEneBan-EL},  with a similar argument as for \eqref{min-rslope-attained}. 
 For
  \eqref{measurable-selection},  observe that, as a consequence of
  Hypothesis \ref{hyp:EL-closedness}, the multivalued mapping
\begin{equation}
\label{upsilon-map}
\mathfrak{A}_{\calR}: (0,\infty)\ti \Spx \rightrightarrows  \Spx^*; \quad  
\mathfrak{A}_{\calR} (\sigma, u):=  \mafo{argmin} 
  \bigset{\calR^*({-}\xi)}{\xi\in \EL\sigma{u}}
\end{equation}
is upper semicontinuous with respect to convergence in $\R$ for $\sigma$ and
weak convergence for $u$.  Therefore, for every measurable selection
$(0,\infty) \ni \sigma \mapsto \wt u_\sigma \in J_\sigma $ the
multivalued mapping $\sigma \mapsto \mathfrak{A}_{\calR}(\sigma, \wt u_\sigma)$
is measurable, as it is given by the composition of measurable mappings. Then,
\cite[Thm.\,3.22]{Castaing-Valadier77} grants the existence of a
\emph{measurable} selection
$\sigma\mapsto \wt\xi_\sigma \in \mathfrak{A}_{\calR} (\sigma, \wt u_\sigma)$.
\end{proof}

We conclude this section with an example in which $\calR$ is not radially
differentiable but  the  De Giorgi identity still holds.

\begin{example}
\label{ex:4}
\upshape
We consider the  $\GBGS$  $(X,\calE,\calR)$ with 
\[
X=\R, \quad \calE(u) = \frac12\, u^2, \quad \calR(v) = \begin{cases} 
\frac12\, v^2&\text{for } |v|\leq 1, \\ 2v^2-\frac32&\text{for }|v|\geq 1.
\end{cases}
\]
Clearly, $\calR$ is not radially differentiable, hence \eqref{eq:BGS.StrangeCond}
does not hold.

For $\ugi=6$ we obtain the unique variational interpolant
\[
\wt u_\sigma= \begin{cases} \frac{24}{4+\sigma} & \text{for } \sigma\in [0,2],\\
6-\sigma& \text{for }\sigma\in [2,5],\\ \frac6{1+\sigma}& \text{for } \sigma\geq 5.
\end{cases}
\] 
It can easily be checked that    identity \eqref{eq:BGS.DeGiEstim} holds for all
$\sigma>0$, where the choice $\wt\xi_\sigma = \rmD\calE(\wt u_\sigma)=\wt u_\sigma$ is
mandatory because the Fr\'echet subdifferential of $\calE$ is single-valued. 
For $\sigma \in [2,5]$ we have $-\wt\xi(\sigma)= -6{+}\sigma    \in
\pl\calR(-1)=[-4,-1]$, where we used $\frac1\sigma(\wt u_\sigma{-}6)=-1$. 

 Recalling the definitions $\delta^\pm_\sigma$ from
\eqref{eq:BGS.MargiEst} we obtain the strict inequalities
$\delta^-_\sigma =-\calR^*(-4)=- \frac 7 2 \lneqq
\partial_\sigma\phi(\sigma)=\sigma{-} \frac{11} 2 \lneqq - \frac 1 2 
=-\calR^*(-1) = \delta^+_\sigma$ for $\sigma \in (2,5)$. 
\end{example}

\subsection{Regularity of  the variational interpolant}
\label{ss:regularity-added}

In this section we take a slight detour from the main theme of the paper and
provide some sufficient conditions for gaining extra time regularity of the
variational interpolant $\sigma \mapsto \wt u_\sigma$. This deviation will only
be useful if we employ a slightly strengthened version of radial
differentiability. 

Furthermore, on the one hand, we will require uniform convexity of the mapping
$u \mapsto \Phi_\sigma (u) $, a sufficient condition for which is, of
course, uniform convexity of the energy functional $\calE$. On the other hand,
we will need to reinforce radial differentiability. Indeed, while the latter
property is equivalent to the fact that the composed, a priori multivalued,
mapping $\calR^*{\circ}\partial\calR : \Spx \rightrightarrows \Spx^*$ is
single-valued (cf.\ Proposition \ref{prop:equivalence-structural}), we will now
further require that $\calR^*{\circ}\partial\calR $ is Lipschitz continuous (on
bounded subsets of $\Spx$).

We collect these conditions in the following
\begin{hypothesis}
\label{h:4reg}
We assume that 
\begin{enumerate}
\item
 There exist $\overline{\lambda}>0$ and $\sigma_*>0$ such that  for all $\sigma>0$
    the mapping $u \mapsto  \Phi_\sigma (u) $
is $\overline\lambda$-convex, namely
\begin{equation}
\label{lambda-cvx-Phisigma}
\begin{aligned}
  & \hspace{-0,5cm} \exists\, \overline\lambda > 0 \ \ \forall\, \sigma >0 \,, \
  \forall\, u_0,u_1\in \Spx \,, \ \forall\, \theta \in [0,1]\,: \\ &
  \hspace{-0,5cm} \Phi_\sigma ((1{-}\theta)u_0 {+} \theta u_1) \leq
  (1{-}\theta)\Phi_\sigma(u_0) + \theta \Phi_\sigma(u_1) -
  \frac{\overline\lambda}2 \theta (1{-}\theta) \|u_0{-}u_1\|^2\,.
\end{aligned}
\end{equation}
\item The mapping $\calR^*{\circ}\partial\calR : \Spx \rightrightarrows
  \Spx^*$  is  Lipschitz continuous on bounded sets, i.e.,
\begin{equation}
\label{geom-cond}
\begin{aligned}
  & \hspace{-0,5cm} \forall\, M>0 \ \exists\, C_M>0 \ \forall\, (v_1,\xi_1), \,
  (v_2,\xi_2) \in \Spx{\times}\Spx^*, \ \max_{i=1,2} \|v_i\|\leq M , \ \xi_i
  \in \partial \calR(v_i) \,:
  \\
  & \qquad \qquad \qquad |\calR^*(\xi_1){-}\calR^*(\xi_2) | \leq C_M \|
  u_1{-}u_2\|\,.
\end{aligned}
\end{equation}
\end{enumerate}
%
%
\end{hypothesis}

Before stating our result, let us pin down two key consequences of Hypothesis
\ref{h:4reg}:
\begin{itemize}
\item It follows from \eqref{lambda-cvx-Phisigma} (cf., e.g., \cite[Prop.\
  2.4]{MiRoSa13NADN}) that the Fr\'echet subdifferential $\partial \Phi_\sigma
  (\cdot): \Spx \rightrightarrows
  \Spx^*$ can be characterized in terms of the following \emph{global}
  estimate: for all $ \sigma>0$ and for every $u \in \dom( \partial
  \Phi_\sigma (\cdot))$ we have that
\begin{equation}
\label{global-characterization}
\begin{gathered}
  \omega \in \partial \Phi_\sigma (u) \text{ if and only if } \\
  \Phi_\sigma (v)- \Phi_\sigma(u) \geq \langle \omega, v{-}u \rangle
  + \frac{\overline\lambda}2 \|v{-}u \|^2 \quad\text{for all } v \in \Spx\,.
\end{gathered}
\end{equation}
\item By Prop.\ \ref{prop:equivalence-structural}, \eqref{geom-cond} in
  particular implies that for every $M>0$ the restriction of
  $\calR$ to the ball
  $\overline{B}_M(0)$ is radially differentiable. Therefore, recalling Lemma
  \ref{l:radial-derivatives} we have that for every $v\in
  \overline{B}_M(0)$ the function
\begin{equation}
\label{deriv-g}
\begin{gathered}
  g(v;\rho) = \rho \calR \left( \frac1\rho v \right) \text{ is differentiable
    at every } \rho>0, \text{ with }
  \\
  \frac{\dd}{\dd \rho} g(v;\rho) = \calR \left( \frac1\rho v \right) -\langle
  \eta, \frac1\rho v \rangle = \calR^*(\eta) \qquad \text{for all } \eta \in
  \partial \calR\left( \frac1\rho v \right)\,.
\end{gathered}
\end{equation}
\end{itemize}

We are now in a position to state the main result of this section.

\begin{theorem}
  Let the $\GBGS$ $(X,\calE,\calR)$ satisfy Hypotheses \ref{h:X} and
  \ref{h:4reg}; and let $\ugi\in \mafo{dom}(\calE)$ be fixed.  Then, for all
  measurable selection
  $(0,\infty) \ni \sigma \mapsto \wt u_\sigma\in J_\sigma $ we have
  $\wt u \in \mathrm{Lip}_{\mathrm{loc}}(]0,\infty[;\Spx)$.
\end{theorem}
\begin{proof}
Preliminarily, we observe that, because of $\sup_{\sigma>0} \sigma \calR((\wt
u_\sigma{-}\ugi)/\sigma) \leq C$, we have
\begin{equation}
\label{curve-bounded}
\exists\, M>0 \ \forall\,  \sigma \in (0,\infty)\, : \quad \| \wt u_\sigma\|\leq M.
\end{equation}
Let us fix $\sigma_{\#}>0$ and let
$[\sigma_1,\sigma_2 ] \subset [\sigma_{\#},\infty)$ be arbitrary.  We apply
\eqref{global-characterization} with $\sigma=\sigma_1$, $u= \wt u_{\sigma_1}$,
$v = \wt u_{\sigma_2}$ and $\omega =0$ as, indeed,
$0 \in \partial \Phi_{\sigma_1} ( \wt u_{\sigma_1})$ is the Euler-Lagrange
equation for $\wt u_{\sigma_1} \in J_{\sigma_1}(\ugi)$. Thus, we obtain
\[
\begin{aligned}
  \frac{\overline\lambda}2 \|\wt u_{\sigma_2}{-}\wt u_{\sigma_1} \|^2 & \leq
  \Phi_{\sigma_1}(\wt u_{\sigma_2})-\Phi_{\sigma_1}(\wt u_{\sigma_1})
  \\
  & = \calE (\wt u_{\sigma_2}) - \calE (\wt u_{\sigma_1}) + \sigma_1 \calR
  \left(\frac1{\sigma_1}(\wt u_{\sigma_2}{-} \ugi) \right) -\sigma_1 \calR
  \left(\frac1{\sigma_1}(\wt u_{\sigma_1}{-} \ugi) \right)\,.
 \end{aligned}
\]
Interchanging $\sigma_1$ and $\sigma_2$ and adding the inequalities gives
\begin{align}
  &\nonumber
  \overline\lambda \|\wt u_{\sigma_2}{-}\wt u_{\sigma_1} \|^2 
  \\
  &\nonumber
  \leq 
  \sigma_2 \calR \big(\frac1{\sigma_2}(\wt u_{\sigma_1}{-} \ugi)  \big)
  -\sigma_1 \calR \big(\frac1{\sigma_1}(\wt u_{\sigma_1}{-} \ugi)  \big) 
  +
  \sigma_1 \calR \big(\frac1{\sigma_1}(\wt u_{\sigma_2}{-} \ugi)  \big)
  -\sigma_2 \calR \big(\frac1{\sigma_2}(\wt u_{\sigma_2}{-} \ugi)  \big) 
  \\
  &
  \nonumber
  =\int_{\sigma_1}^{\sigma_2}  \left[ \frac{\dd}{\dd \rho} g(\wt
    u_{\sigma_1};\rho){-} \frac{\dd}{\dd \rho} g(\wt u_{\sigma_2};\rho)
  \right] \dd \rho 
  \overset{\text{\eqref{deriv-g}}}= 
  \int_{\sigma_1}^{\sigma_2} \left( \calR^*
    (\eta_1^\rho){-}\calR^* (\eta_2^\rho) \right)  \dd \rho
\end{align}
for every $\eta_i^\rho \in \partial \calR (\tfrac1\rho \wt u_{\sigma_i})$,
$i=1,2$. Now using the Lipschitz property, we find
\begin{align}
  & \nonumber
  \overline\lambda \|\wt u_{\sigma_2}{-}\wt u_{\sigma_1} \|^2  \leq
  \int_{\sigma_1}^{\sigma_2} C_M\| \tfrac1\rho \wt u_{\sigma_2}{-}  \tfrac1\rho
  \wt u_{\sigma_1} \|\dd \rho  
\end{align}
due to \eqref{geom-cond}.  All in all, using that $\sigma_1 \geq \sigma_{\#}$ we
conclude that   
\begin{equation}
  \label{key-est-4AC}
  \overline\lambda \|\wt u_{\sigma_2}{-}\wt u_{\sigma_1} \|^2  \leq
  \frac{C_M}{\sigma_{\#}} (\sigma_2{-}\sigma_1) \|\wt u_{\sigma_2}{-}\wt
  u_{\sigma_1} \| \, ,
\end{equation}
hence $\sigma \mapsto \wt u_{\sigma} $ is Lipschitz continuous on
$[\sigma_{\#},\infty)$. By the arbitrariness of $\sigma_{\#}$, we conclude
$\wt u \in \mathrm{Lip}_{\mathrm{loc}}(]0,\infty[;\Spx)$.
\end{proof}

\subsubsection*{More on condition \eqref{geom-cond}.} 
We conclude this section by gaining further insight into property
\eqref{geom-cond}. The key observation is that it  is implied by  radial
differentiability \text{and} superlinearity of $\calR^*$.

\begin{lemma}
\label{l:characteriz-geom2}
Let $\calR: X \to [0,\infty)$ satisfy Hypothesis \ref{h:X},  be
radially differentiable, \emph{and} suppose that $\calR^*$ has superlinear growth at infinity, i,e.,
\begin{equation}
\label{super-Diss-star}
\lim_{\|\xi\|_*\to\infty} \frac{\calR^*(\xi)}{\|\xi\|_*} = \infty\,.
\end{equation}
Then, $\calR^*{\circ}\partial\calR : \Spx \to \Spx^*$
fulfills \eqref{geom-cond}.
\end{lemma}
\begin{proof}
 It follows from  
 \eqref{super-Diss-star} that the subdifferential $\partial\calR: X \rightrightarrows X^*$ is a bounded operator,  cf.\   \cite[Prop.\ 2.14]{Brez73}. 
Then,  \eqref{geom-cond} immediately follows, taking into account that $\calR^*$ is itself  Lipschitz on bounded sets since $\calR$  has superlinear growth at infinity.
\end{proof} 
We also remark that, if $\calR^*$
 is superlinear, then property \eqref{geom-cond} is equivalent to the  Lipschitz continuity on bounded sets of the mapping $\Prod$ from \eqref{equal-result}. This follows from the fact that 
 $\Prod(u) = \calR(u) + (\calR^*{\circ}\partial\calR)(u)$, and $\calR$ is  Lipschitz continuous on bounded sets.  Exploiting this observation, we can readily prove that 
property
\eqref{geom-cond} is stable under the sum of dissipation potentials. 

\begin{corollary}
\label{cor:stable-under-sum}
Let $\calR_i: X \to [0,\infty)$, $i=1,2$, fulfill Hypothesis \ref{h:X} be
radially differentiable,  and satisfy \eqref{super-Diss-star}. Then, $\calR_1+\calR_2$ satisfies \eqref{geom-cond}.
\end{corollary}
\begin{proof}
  It suffices to remark (with slight abuse of notation) that
\[
  \Prod_{\calR_1{+}\calR_2}(v) = \langle
  \partial(\calR_1{+}\calR_2)(v),v\rangle = \langle \partial \calR_1(v)
  {+}\partial \calR_2 (v),v\rangle = \Prod_{\calR_1}(v)+\Prod_{\calR_2}(v)
  \quad \text{for all } v \in X\,,
\]
where the validity of the sum rule
$\partial(\calR_1{+}\calR_2)(v) = \partial \calR_1(v) {+}\partial \calR_2 (v)$
is guaranteed, under the present assumptions, by, e.g., \cite[Thm.\,1,
p.\,211]{IofTih79TEPe}.
\end{proof}

Ultimately, $p$-homogeneous dissipation potentials  provide
examples for the validity of~\eqref{geom-cond}.

\begin{corollary}
\label{co:pHomog}
  Let $\calR: \Spx \to [0,\infty)$ be positively homogeneous of degree $p$ for
  some $p> 1$. Then, 
  \begin{equation}
\label{Euler-relation}
\Prod(v) = p \calR(v) \qquad \text{for all } v \in \Spx
\end{equation}
and 
  $\calR^*{\circ}\partial\calR$ fulfills \eqref{geom-cond}.
\end{corollary}
\begin{proof}
  First of all, we show \eqref{Euler-relation} (which, in fact, also holds for $p=1$).   Indeed, 
  it suffices to observe that, by definition of
$\partial \calR(v)$ and $p$-homogeneity of $\calR$, we have for all
$\xi \in \partial \calR(v)$
\[
  (\lambda^p {-}1)\calR(v) = \calR(\lambda v) - \calR(v)\geq \langle \xi,
  \lambda v{-}v \rangle = (\lambda{-}1) \Prod(v)\,.
\]
Then, we divide the above estimate by $ (\lambda{-}1) $ for all $\lambda>1$,
and take the limit as $\lambda \to 1^+$, thus obtaining
$p \calR(v) \geq \Prod (v)$. The converse inequality follows by dividing by
$ (\lambda{-}1) $ for all $0<\lambda<1$ and sending $\lambda \to 1^-$.

From \eqref{Euler-relation} and the  fact that $\calR$ is Lipschitz
continuous on bounded sets, we deduce that $\Prod$ satisfies
\eqref{geom-cond}. 
\end{proof}

Eventually, by Corollary \ref{cor:stable-under-sum} property \eqref{geom-cond}
holds for linear combinations of homogeneous potentials of possibly different
degrees, as well.

\begin{example}[Example \ref{PDE-example} re-revisited] 
\slshape  Getting back to Example \ref{PDE-example}, let us additionally 
assume that $\mathrm{R}$ is differentiable on $\R {\setminus} \{0\}$: then 
the dissipation potential $\calR(v) = \int_\Omega \mathrm{R}(v(x)) \dd x$ is 
radially differentiable by Remark \ref{rmk:ex-revisited}.  Since we have required
$\mathrm{R}^*$ to have $p'$-growth with $p'>1$, $\calR^*$ has superlinear
growth at infinity, and then Lemma \ref{l:characteriz-geom2} applies, ensuring
that property \eqref{geom-cond} holds. 
\end{example}

\subsection{De Giorgi estimates for general  dissipation potentials}
\label{ss:4.2}

We now return to the general case in which $\calR$ is allowed to take the value
$\infty$ with an open domain $\dom(\calR)$ in $\Spx$ (cf.\ Hypothesis
\ref{h:X}).  We will show that, if we drop the radial differentiability
condition \eqref{eq:BGS.StrangeCond} on $\calR$,  the  De Giorgi simple
and improved estimates can be still retrieved for $\GBGS$ $\BGS$ featuring
general dissipation potentials $\calR$, see Theorem
\ref{th:DGL.NON.RD.Simple} ahead.

Our strategy for proving both results  consists in replacing $\calR$ 
by the Yosida approximations $(\calR_\eta)_\eta$,  where we use a suitable
norm on the reflexive Banach space $\Spx$ such that all $\calR_\eta$ are
differentiable. Thus, for all $\eta>0$ we find measurable variational
interpolants 
\begin{equation}
\label{min-eta}
\sigma \mapsto \wt u^\eta_\sigma \in \Jugi_\eta(\sigma):=
\mafo{Argmin}\Bigset{\sigma\calR_\eta\big( \tdfrac1\sigma(u{-}\ugi)\big) +
  \calE(u)}{ u \in X } 
\end{equation}
satisfying  the  De Giorgi identity for the corresponding $\GBGS$
$(X,\calE, \calR_\eta)$. We will then study the limit $\eta \to 0^+$.
      
A related approach based on Yosida regularization was developed in
\cite[Chap.\,3]{Bacho:PHD} in order to extend the existence result for $\GBGS$
provided in \cite{MiRoSa13NADN}, \emph{without} the radial differentiability
condition that was adopted therein (recall the discussion prior to Prop.\
\ref{prop:equivalence-structural}).  In \cite{Bacho:PHD} the limit
$\eta\to 0^+$ was performed on the level of the solutions of the evolution
equation.  

Here, we will resort to Yosida approximation to get rid of radial
differentiability for the very proof of  the  De Giorgi estimate.
Hence, the limiting process $\eta \to 0^+$ will be more delicate than in
\cite{Bacho:PHD}, because variational interpolants are not absolutely
continuous in $[0,\sigma_0]$, in contrast to solutions of the regularized
$\GBGS$.   Nevertheless, starting from our approximations
$(\sigma\mapsto \wt u^\eta_\sigma)_{\eta>0}$ we will be able to find a limiting
function $\sigma \mapsto \wt u_\sigma \in \Spx $ that such that
\[
\forall \, \sigma\in (0,\sigma_*)\ \exists\, (\eta^\sigma_k)_{k\in \N}: \quad 
\eta^\sigma_k \to 0 \ \text{ and } \ \wt u^{\eta^\sigma_k}_\sigma \weak \wt
u_\sigma \quad \text{for }k\to \infty,
\]
see Proposition \ref{pr:MeasSelect}.  This pointwise weak convergence along a
$\sigma$-dependent subsequence will the crucial point to pass to the liminf, in
 the  De Giorgi identity for $(X,\calE, \calR_\eta)$, to derive the De
Giorgi estimates,  also based on the \textsc{Mosco}-convergence as
$\eta \to 0^+$ of the functionals
\begin{equation}
\label{Phi-eta}
\Phi_\sigma^\eta (u): =\sigma\calR_\eta\big( \tdfrac1\sigma(u{-}\ugi)\big) +
  \calE(u) 
  \end{equation}
    involved in \eqref{min-eta}, to the functional $\Phi_\sigma$ for fixed $\sigma>0$, which means 
\begin{equation}
\label{Mosco-Phi-eta}
\left\{
\begin{array}{ll}
\forall\, (u_\eta)_\eta, u  \in  X \,:  \ u_\eta \rightharpoonup u \ \Rightarrow \ \liminf_{\eta \to 0^+} \Phi_\sigma^\eta (u_\eta) \geq \Phi_\sigma(u),
\\
\forall\, u \in X \ \exists\, (\hat{u}_\eta)_\eta\, : \hat{u}_\eta \to  u \text{ and }  \limsup_{\eta \to 0^+} \Phi_\sigma^\eta (\hat{u}_\eta) \leq \Phi_\sigma(u).
\end{array}
\right.
\end{equation}
This will follow easily from the Mosco-convergence of the potentials
$(\calR_\eta)_\eta$ stated in the upcoming Lemma \ref{le:Yosida.calR}.  
      
For the proof of the simple De Giorgi estimate it will be sufficient to invoke
the lower semicontinuity of the $\calR$-slope of $\calE$, namely $\rsloname$,
 in the limit passage in  the  De Giorgi identity holding for the
regularized potentials.  To obtain the improved  the  De Giorgi
estimate, we will have to impose an additional condition on the lower
semicontinuity of the conditioned $\calR$-slope, namely $\wrsloname$:

\begin{hypothesis}[Qualified lower semicontinuity of $\wrsloname$]
\label{hy:QualLSC.CR}
Fix $\sigma \in (0,\sigma_*)$ and $ \ugi \in \dom(\calE)$ and consider $u_\eta \in
\Jugi_\eta(\sigma)$ for $\eta \in {]0,1[}$. Then,
\begin{equation}
  \label{eq:LSC.CR}
  u_\eta \weak u \quad \Longrightarrow \quad \liminf_{\eta\to 0}
  \wrslopar{\eta}{u_\eta}\sigma \geq  \wrslo{u}\sigma .
\end{equation}
\end{hypothesis}

We will discuss sufficient conditions for this hypothesis after stating the
main result of this subsection. We emphasize that the following result is
weaker than  the  De Giorgi identity in two ways. First we only have
the inequalities instead of the equality, and second we only prove the
existence of \emph{at least one} measurable interpolant, whereas  the 
De Giorgi identity holds for \emph{all} measurable interpolants with
$\wt u_\sigma \in \Jugi(\sigma)$.   So far, we are not able to show that
all measurable interpolants have this property.
  
 Let us also mention in advance that, mirroring the statement of Theorem
\ref{th:DGL.RadialDiff}, under the additional Hypothesis
\ref{hyp:EL-closedness}, we are in a position to prove the existence of a
measurable force selection for which the  improved De Giorgi 
estimate holds.  An analogous statement could be given for  the simple
De Giorgi  estimate; we choose to omit it to avoid overburdening the
exposition.

\begin{theorem}[De Giorgi estimates for general $\calR$]
\label{th:DGL.NON.RD.Simple} 
Let the $\GBGS$ $(X,\calE,\calR)$ satisfy Hypothesis \ref{h:X}, and let
$\sigma_*>0$.  Then, for every $\ugi\in \mafo{dom}(\calE)$ there exists a
measurable variational interpolant
$\sigma \mapsto \wt u_\sigma\in J_\sigma \subset X$
fulfilling  the  \emph{simple De Giorgi estimate}  
\begin{align}
  \label{eq:BGS.DeGiIineqSIMPLE}
  \forall\, \sigma \in (0,\sigma_*)\colon \quad \calE(\wt u_\sigma) 
  + \sigma \calR\big( \frac1\sigma(\wt u_\sigma{-}\ugi)\big) 
  + \int_0^\sigma \calS_\calR(\wt u_\sigma) \dd \rho \leq \calE(\ugi).   
\end{align}
f Hypothesis \ref{hy:QualLSC.CR} holds, then this interpolant
also satisfies  the   \emph{improved De Giorgi estimate}%
\begin{subequations}
\begin{align}
  \label{eq:BGS.DeGiIineqImprov}
  \forall\, \sigma \in (0,\sigma_*)\colon \quad \calE(\wt u_\sigma) 
  + \sigma \calR\big( \frac1\sigma(\wt u_\sigma{-}\ugi)\big) 
  + \int_0^\sigma   \wrslo{ \wt u_\sigma}\sigma \dd \rho \leq \calE(\ugi).   
\end{align}
 Thus, under the additional Hypothesis \ref{hyp:EL-closedness} there
exists a measurable force selection
$\sigma \mapsto \wt\xi_\sigma\in \EL {\sigma}{\wt u_\sigma} $ satisfying the
\emph{improved De Giorgi estimate with force selection}
\begin{align}
  \label{eq:BGS.DeGiIneq-wsel}
  \forall \, \sigma>0: \quad \calE(\wt u_\sigma) + \sigma \calR\big( \frac1\sigma(\wt
  u_\sigma{-}\ugi)\big) + \int_0^\sigma \calR^*(-\wt\xi_\rho) \dd \rho \leq
  \calE(\ugi) \, . 
\end{align}
\end{subequations}
\end{theorem} 

We postpone the proof of this main result to Section \ref{su:PROOFS}. 

 Prior to discussing the  rather abstract lower semicontinuity property
\eqref{eq:LSC.CR} on the conditioned slope $\wrsloname$,  we need to pin
down some of the properties of the Moreau-Yosida approximations
$(\calR_\eta)_{\eta>0}$.  First of all, let us settle their definition: 
Since $X$ is reflexive, Asplund's renorming theorem \cite{Asplund68}, see also
\cite[Thm.\,1.105]{BaPre12}, ensures that there exists an equivalent norm on
$X$ (still denoted by $\| \cdot \|_X$), such that, correspondingly, both $X$
and $X^*$ are strictly convex (namely, the mappings $x\mapsto \|x\|_{X}^2$ and
$\xi\mapsto \|\xi\|_{X^*}^2$ are strictly convex). With this, the Moreau-Yosida
approximations $(\calR_\eta)_{\eta>0}$ of $\calR$ are defined via
\[
  \calR_\eta \colon \Spx \to [0,\infty) \qquad \calR_\eta(v): = \inf_{w\in
    \Spx} \left( \frac{\|w{-}v\|_X^2}{2\eta} {+}\calR(w)\right)\,.
\]
The following result collects the properties of the functionals
$(\calR_\eta)_{\eta>0}$ we are going to use.

\begin{lemma}
\label{le:Yosida.calR}
For all $\eta>0$ the functional $\calR_\eta$ is differentiable on $\Spx$ and
has the convex conjugate $\calR_\eta^*(\xi) = \calR^*(\xi) +\frac\eta2 
\|\xi\|^2_{X^*} $ for $ \xi \in \Spx^*$.  Moreover,
\begin{enumerate}
\item the family $\big(\calR_\eta\big)_{\eta\in (0,1]}$ is equi-coercive, i.e.\
  \begin{equation}
    \label{equi-coercivity-Reta}
    \forall\, S>0 \ \ \exists\, C_S>0 \ \ \forall\, \eta\in (0,1] \  \
    \forall\, v \in \Spx : \quad  
    \calR_\eta(v) \leq S  \ \ \Longrightarrow \ \   \|v\|\leq C_S\,;
  \end{equation}
\item the family $(\calR_\eta)_\eta$ \textsc{Mosco}-converge to $\calR$ as
  $\eta \downarrow 0$;  in particular, we have 
\[
  v_\eta \weak v \ \Longrightarrow \ \calR(v) \leq \liminf_{\eta\to 0}
\calR_\eta(v_\eta) \quad \text{and} \quad 
\xi_\eta \weak \xi \ \Longrightarrow \ \calR^*(\xi) \leq 
 \liminf_{\eta\to 0} \calR^*_\eta(\xi_\eta).  
\] 
\item  For all $(v_\eta)_\eta \subset X$,  $(\xi_\eta)_\eta \subset X^*$
\begin{equation}
\label{closedness-eta}
\left. \ba{cl} 
v_n \to v & \text{as }\eta \to 0^+,
\\
\xi_\eta \rightharpoonup \xi & \text{as }\eta \to 0^+, 
\\
\xi_\eta\in \partial\calR_\eta (v_\eta) & \text{for all }\eta>0
\ea \right\} \quad \Longrightarrow \quad  \xi \in \partial \calR(v)\,.
\end{equation} 
\end{enumerate}\vspace{-0.8em}
\end{lemma}
\begin{proof}
By strict convexity of $X^*$, the space $X$ is \emph{smooth} (cf., e.g.,
\cite[Thm.\,1.101]{BaPre12}), which is equivalent to the property that the
duality mapping $\rmJ_X =  \partial\big(\frac12 \| \cdot \|_X^2\big) 
\colon X \rightrightarrows X^*$ is \emph{single-valued}  and one-to-one,
 see \cite[Rmk.\,1.100]{BaPre12}. Hence, the differentiability of the 
functionals $\calR_\eta$  follows, cf.\ \cite[Sec.\,3.4.1]{Att84VCFO};  in
particular, we have 
\begin{equation}
  \label{psi-lambda}
  \calR_\eta(v) = \calR\big(W_\eta(v) \big) +
  \frac{\|W_\eta(v){-}v\|^2}{2\eta} \ \text{ and } \ 
  \rmD \calR_\eta(v) = \frac1\eta \rmJ_X\big(v{-}W_\eta(v)\big) 
  \quad \text{for } v  \in X,
\end{equation}
with $W_\eta(v) = 
\big(\mathrm{Id}{+}\eta\rmJ_X{\circ}\partial\calR\big)^{-1}(v)$. 

To show \eqref{equi-coercivity-Reta}, we use \eqref{psi-lambda} and that
$\calR$ has superlinear growth. This gives
\[
  \exists\, \overline C>0 \ \ \forall\, \eta>0 \ \ \forall\, v \in X: \quad
  \| W_\eta(v) \|\leq \overline{C} + \calR(W_\eta(v)) \leq
  \overline{C} + \calR_\eta(v) \,.
\]
 Moreover, using $\calR\geq 0$ we find  
\[
  \|W_\eta(v){-}v\|  \leq ( 2\eta \calR_\eta(v))^{1/2}
\leq ( 2 \calR_\eta(v))^{1/2} \quad \text{for all }
  \eta \in (0,1]\,.
\]
Thus,
$ \|v \| \leq \| W_\eta(v) \| +\|W_\eta(v){-}v\| \leq \overline C +
\calR_\eta(v) + (2\calR_\eta(v))^{1/2}$, and \eqref{equi-coercivity-Reta}
follows.

The \textsc{Mosco} convergence follows from the results in
\cite[Sec.\ 3.4]{Att84VCFO}.
Finally, 
 in order to show \eqref{closedness-eta},  we note that $\xi_\eta\in \partial\calR_\eta (v_\eta)
$ 
equivalent to 
\[
\big\langle \xi_\eta , v_\eta \big\rangle  
= \calR_{\eta} \big( v_\eta \big) + \calR^*_{\eta}(\xi_\eta) 
= \calR_{\eta_n} \big( v_\eta\big) + \calR^*(\xi_\eta) +
\tfrac{\eta}2 \| \xi_\eta\|_{X^*}^2 .  
\]
Using that $v_\eta \to v$ and $\xi_\eta\weak  \xi$  we can pass to the
limit on the 
left-hand side. On the right-hand side we can pass to the liminf and obtain
$\langle \xi , v\rangle \geq \calR(v)+\calR^*(\xi) $, which yields 
 the
desired result $\xi\in \pl \calR(v)$.  
\end{proof}

 The following result 
shows that
 Hypothesis \ref{hy:QualLSC.CR} can be naturally deduced if  the conditioned  subdifferential $\ELname
\calE$ has a suitable closedness property.
A  sufficient condition for such closedness  involves the following property, which we define in general: 
 We say that a functional
$\calF:X\to (-\infty,\infty]$ has the  \emph{weak-implies-strong property}
($\mathrm{WIS}$, for short), 
 if the following
implication holds:
\begin{equation}
\label{eq:WeakStrongP}
u_k\overset{X}\weak  u \text{ (weakly) } \text{ and } \calF(u_k)\to \calF(u)<\infty \quad
\Longrightarrow \quad u_k \overset{X}\to u \text{ (strongly)}.
\end{equation}
This property is certainly true if $\calF$ has compact sublevels, however it is
also true for uniformly convex functions, or for functionals $\calF(u)= g(\|
u\|)$ if $g$ is strictly increasing and $\|\cdot\|$ is a uniformly convex norm.   Thus, the conditions in Prop.\ \ref{pr:SuffC.CR.lsc} are sufficiently 
easy to check, and indicate the potential for generalizations.

\begin{proposition}[Sufficient conditions for \eqref{eq:LSC.CR}]
\label{pr:SuffC.CR.lsc}
Let the $\GBGS$ $(X,\calE,\calR)$ satisfy Hypothesis \ref{h:X} and fix
$\sigma_*>0$ and $\ugi \in \dom(\calE)$. 

(A) If for all $\sigma\in (0,\sigma_*)$, the  operator $\ELname \calE:  \mathrm{Graph}(\Jugi)  \rightrightarrows \Spx^*$ satisfies 
\begin{align}
\label{eq:HypEneBan-EL} 
{\left. \begin{array}{c}  (\eta_n,  u_n,\xi_n) \rightharpoonup
 (0,   u,\xi) \ \text{ in }  [0,\infty)  \ti \Spx \ti
   \Spx^*, \\  u_n \in    \Jugi_{\eta_n} (\sigma)  
\text{ for all } n \in \N, \\ \text{ and } \xi_n \in
 \pl_{\calR_{\eta_n}}  \calE(\sigma; u_n) \ \text{ for all } n \in \N
\end{array} \right\} \ \Longrightarrow \ \xi \in \EL{\sigma}{u}\,,}
\end{align}
 then the lower semi-continuity property \eqref{eq:LSC.CR} is satisfied. 

(B) If $\calE$ or $\calR$ have the   $\mathrm{WIS}$   property
and $\mathrm{dom}(\calR)=X$ with $\calR^*$ having superlinear growth at infinity, 
 then
condition \eqref{eq:HypEneBan-EL} holds.
\end{proposition}
\begin{proof} 
For (A) the proof is straightforward.  We note that for a sequence $(\eta_n,u_n)\weak (0,u)$ with 
$\alpha:= \liminf_{n\to \infty}   \wrslopar{\eta_n}{u_\eta}\sigma < \infty$ we
find  $\xi_n\in   \pl_{\calR_{\eta_n}}  \calE(\sigma; u_n) \subset 
X^*$ with  $ \calR_{\eta_n}^*(\xi_n) \leq \wrslopar{\eta_n}{u_\eta}\sigma + \tdfrac 1 n$. In turn,
since $ \calR^*(\xi_n) \leq 
 \calR_{\eta_n}^*(\xi_n)
$, by the coercivity property \eqref{coercivityRstar} we conclude that 
$  \|\xi_n\|_{X^*} \leq C< \infty $.

Thus, after extracting a suitable
subsequence (not relabeled) we may assume $\xi_n\weak \xi$ and the closedness
property \eqref{eq:HypEneBan-EL} implies $\xi \in
\EL{\sigma}{u}$. Now, using the weak  lower semicontinuity of $\calR^*$
we have 
\begin{align*} 
\wrslopar{}{u}\sigma & \leq \calR^*(\xi) \leq \liminf_{n\to \infty}
\calR^*_{\eta_n}(\xi_n) = \liminf_{n\to \infty} \big(
\wrslopar{\eta_n}{u_\eta}\sigma +\tdfrac1n \big) = \alpha,
\end{align*}
which is the desired result of part (A). 

For part (B),   let $(\eta_n,  u_n,\xi_n)$ be as in  \eqref{eq:HypEneBan-EL}.  
We use that $u_n \in \Jugi_{\eta_n}(\sigma)$ and $u_n
\weak u$ already imply  
\begin{equation}
\label{convergence-of-sum}
\Phi_\sigma^{\eta_n}(u_n) = \sigma\calR_{\eta_n}\big( \tdfrac1\sigma(u_n{-}\ugi)\big) +
  \calE(u_{\eta_n}) 
    \longrightarrow  \Phi_\sigma(u)\,.
\end{equation}
  Indeed, the $\liminf$-estimate follows by 
\eqref{Mosco-Phi-eta}, whereas for the $\limsup$ we observe that 
$\Phi_\sigma^{\eta_n}(u_n) \leq \Phi_\sigma^{\eta_n}(u) \leq  \Phi_\sigma(u)$, where the first inequality follows from $u_{\eta_n} \in  \Jugi_{\eta_n}(\sigma)$ and the second one from $\calR_{\eta_n} \leq \calR$. 
We also have the individual $\liminf$-estimates 
\begin{equation}
\label{individual-liminfs}
\calR(\frac1\sigma(u{-}\ugi)) \leq \liminf_{n \to   \infty} \calR_{\eta_n} 
  (\frac1\sigma(u_n{-}\ugi))  \quad \text{and} \quad 
\calE(u) \leq \liminf_{n \to  \infty} \calE(u_n),
\end{equation}
the first one  by the Mosco-convergence of $\calR_{\eta_n}$ to $\calR$, and the second one by the weak lower semicontinuity of $\calE$. 
Combining \eqref{convergence-of-sum} and \eqref{individual-liminfs}
we see that both liminfs are limits and that the inequalities must be
equalities. Hence, by assumed $\mathrm{WIS}$ property for either $\calR$, or
$\calE$, we find the strong convergence $u_n\to u$.   Now, we need to take
pass to the limit, as $n\to\infty$, in both relations:
$\xi_n\in \pl\calE(u_n) $ and
${-}\xi_n \in \partial\calR_{\eta_n}\big( \frac1\sigma(u_n{-}\ugi)\big)$. For
the first relation, we resort to the closedness property \eqref{HypEneBan-2}
for $\partial \calE$ (observe that $u_n$ lies in an energy sublevel since
$\calE(u_n) \leq \Phi_\sigma^{\eta_n}(u_n) \leq \calE(\ugi)$).  For the second
relation, we exploit the gained strong convergence
$ \frac1\sigma(u_n{-}\ugi) \to \frac1\sigma(u{-}\ugi)$ and property
\eqref{closedness-eta}.  All in all, we conclude $\xi \in \EL{\sigma}{u}$,
and \eqref{eq:HypEneBan-EL} is established.  This proves part (B).
\end{proof}

\subsection{Proof of Theorem \ref{th:DGL.NON.RD.Simple}}
\label{su:PROOFS}

 The first tool in the proof of Theorem \ref{th:DGL.NON.RD.Simple} are the
previously recapped properties of the Moreau-Yosida approximations
$\calR_\eta$.  The second tool is the following standard result on the
existence of measurable selections for Kuratowski upper limits  which,
like Lemma \ref{l:conseq-of-closedness}, also draws on
\cite{Castaing-Valadier77}.  

\begin{proposition}[Measurable selection]
\label{pr:MeasSelect}
Consider a sequence of $(u^k)_{k\in \N}$ of measurable functions
$(0,\sigma_*) \ni \sigma \mapsto u^k_\sigma \in M$, where $M$ is a compact
metric space.
\par Then, there exists a measurable selection
$\sigma \mapsto \wt{u}_\sigma$
of the set of all limit points of $ (u^k_\sigma)_k$, namely
   for all $\sigma\in (0,\sigma_*)$  there
exists a subsequence $(k^\sigma_n)_{n\in \N}$ with $k^\sigma_n \to 
\infty$ and $u^{k^\sigma_n}_\sigma \to \wt{u}_\sigma$ in $M$.

 Furthermore,  if $\calA:(0,\sigma_*)\ti M \to [0,\infty]$ is a measurable
function such that $\calA(\sigma,\cdot)$ is lower semicontinuous for all
$\sigma$, then the measurable selection $\wt u$ can be chosen
such that additionally 
\[
  \calA(\sigma, \wt u_\sigma) \leq  A(\sigma):=\liminf_{k\to \infty}
  \calA(\sigma, u^k_\sigma) \quad \text{for all } (0,\sigma_*).
\]
\end{proposition} 
\begin{proof}  As $M$ is compact, we can define the multivalued
  mapping 
\begin{equation}
  \label{multivalued-U}
  \mathscr{U}:  (0,\sigma_*)\rightrightarrows M, \qquad 
  \mlt{U}\sigma:= \Ls_{k\to \infty} \{ u^k_\sigma \}, 
\end{equation}
where $\Ls$ is the set of all limit points (in the sense of a Kuratowski upper
limit,  cf.\ \eqref{multivalued-usc-2}) 
\begin{equation}
  \label{charact-LS}
w\in \Ls_{k\to \infty} \{u^k_\sigma\} \quad \Longleftrightarrow \quad
\exists\, (k_n)_{n\in \N}: \ k_n\to \infty \text{ and } u^{k_n}_\sigma \to w.  
\end{equation}
Hence, each $ \mlt{U}\sigma$ is well defined, non-empty and closed. The results
of \cite[Thm.\,8.2.5]{Aubin-Frankowska} imply that $\mathscr{U}$ is a
measurable multivalued mapping.  Now, \cite[Thm.\,3.22]{Castaing-Valadier77}
grants the existence of a measurable selection for $ \mathscr{U}$, i.e.\ a
function $\sigma \mapsto \wt u_\sigma$ with $\wt u_\sigma \in \mlt{U}\sigma$
for all $\sigma\in (0,\sigma_*)$.

For the  second part of the statement  involving $\calA$ we define
\[
  \mathscr U_\calA(\sigma):=\set{u \in \mlt{U}\sigma }{ \calA(\sigma,u)\leq
  A(\sigma)}.
\]
By the definition of $\mathscr U_\calA$ and the lsc of $\calA$ we know have
$ \mathscr U_\calA(\sigma)\neq \emptyset$ for all $\sigma $. Hence, 
$\mathscr{U}_\calA: (0,\sigma_*)\rightrightarrows M$ is a measurable
multivalued mapping with non-empty and closed values, and we can apply the
selection theorem as above, see also \cite[Thm.\,8.1.4]{Aubin-Frankowska}.
\end{proof}

We are now ready to carry out  the proof of Theorem \ref{th:DGL.NON.RD.Simple}
which provides  the  De Giorgi simple and improved estimate for general dissipation
potentials $\calR$. For this, consider variational integrands for the
regularized $\GBGS$ with dissipation potential $\calR_\eta$.\medskip

\noindent
\begin{proof}[Proof of Theorem \ref{th:DGL.NON.RD.Simple}] 
\mbox{} \\[0.2em]
\STEP{1: Regularization.} For $\eta\in (0,1]$ we consider the dissipation
potentials $\calR_\eta$, which are radially differentiable. Hence, by Theorem
\ref{th:DGL.RadialDiff} there exist variational interpolants $\sigma \mapsto 
\wt u^\eta_\sigma$ such that  the  De Giorgi identity holds
\begin{align}
  \label{eq:DeGiIdent-eta-slope}
  \forall\,   \sigma>0\,: 
  \quad  \calE(\tyos u\sigma\eta) + \sigma \calR_\eta\big( \frac1\sigma(\tyos
  u\sigma\eta {-}\ugi)\big) + \int_0^\sigma \wrslopar{\eta}{\tyos u\rho\eta}
  {\rho} \dd \rho =  \calE(\ugi)\,.
\end{align} 
In particular, we have $ \tyos u\sigma\eta \in \Jugi_\eta(\sigma)$,  and
$\wrslopar{\eta_n}{\cdot}{\cdot}$ is defined in \eqref{eq:ConditSlope}.\smallskip  

\noindent
\STEP{2: A priori bounds.}  Clearly, we have that
$\calE(\tyos u\sigma\eta)\leq \calE(\ugi)$ for all $\sigma$ and $\eta$.  We
will now use the place-holder
$\yos v\sigma\eta: = \frac1\sigma (u_{\sigma,\eta}{-}\ugi)$.  Recalling that
$\calE$ is bounded below by some constant $ E_0$, we deduce from
\eqref{eq:DeGiIdent-eta-slope} that
\begin{equation}
  \label{eqanother-estimate}
  \sigma\calR_\eta (\yos v\sigma\eta) 
  +E_0 \leq \sigma\calR_\eta (\yos v\sigma\eta) 
  +   \calE(\tyos u\sigma\eta)  \stackrel{(1)}{\leq}   \calE(\ugi) 
  \quad \text{for all } \sigma,\, 
  \eta>0\,,
\end{equation}
where {\footnotesize (1)} derives from the minimality of $\tyos u\sigma\eta$.
Thus, it follows from \eqref{equi-coercivity-Reta} applied with
$S := \calE(\ugi)$ that there exists a constant $C_S>0$ such that
$ \sigma \|\yos v\sigma\eta \| \leq C_S$.  Now, recalling that
$\yos v\sigma\eta= \frac1\sigma (\tyos u\sigma\eta{-}\ugi)$ we infer that
$  \|\tyos u\sigma\eta \| \leq \sigma \|\yos v\sigma\eta\|\ +\| \ugi \| \leq C_S
  + \| \ugi \| $.
All in all, we proved that 
\begin{equation}
  \label{eq:thanks-Alex}
  \exists\,\overline E>0 \ \ \forall\, \eta \in (0,1]\, \ \forall\, 
   \sigma \in (0,\sigma_*)\, :  \qquad   
   \calE(\tyos u\sigma\eta)+ \|\tyos u\sigma\eta \| \leq \overline E.\smallskip 
\end{equation}

\noindent 
\STEP{3: Measurable selection.} Using the result of Step
2, we are able to apply Proposition \ref{pr:MeasSelect} with the choices
$M=B_{\ol E}(0)\cap S_{\ol E} \subset X$ equipped with the weak topology and 
$\calA=\rsloname$ which is lsc, see \eqref{lsc-slope}. Hence, we obtain a 
measurable function $\sigma \mapsto \wt u_\sigma \in X$ with
$ \tyos u\sigma\eta \weak \wt u_\sigma$ along a subsequence
$\eta= \eta^\sigma_k\searrow 0$, where $(\eta^\sigma_k)_{k\in \N}$ may depend
on $\sigma$. Moreover, 
\begin{equation}
  \label{eq:SR.liminf}
  \rslo{\wt u_\sigma} \leq \liminf_{\eta\to 0} \rslo{\tyos u\sigma\eta}
   \quad \text{for all } \sigma\in (0,\sigma_*). 
\end{equation}

\noindent 
\STEP{4: Convergence of $\Phi_\sigma$.}  Recall that 
$\Phi^\eta_\sigma(u)=\calE(u) + \sigma
\calR_\eta\big(\frac1\sigma(u{-}\ugi)\big)$ and that the family
$\big( \Phi^\eta_\sigma(\cdot) \big)_\eta$ Mosco-converges to
$\Phi_\sigma(\cdot)$.  Moreover, the functions
$\sigma \mapsto \phi^\eta(\sigma) = \min\bigset{\Phi_\sigma^\eta(u)}{u\in
  X}$ are given independently of our curves $\tyos u\eta\sigma$ and
$\wt u_\sigma$. 

From $ \tyos u\eta\sigma \in \Jugi_\eta(\sigma)$ we
have $\phi^\eta(\sigma)= \Phi^\eta_\sigma(\tyos u\sigma\eta)$. Using the
convergence $ \tyos u\sigma{\eta^\sigma_k} \weak \wt u_\sigma $ as
$k \to \infty$ and the Mosco-convergence of
$\Phi^\eta_\sigma $ we obtain  (cf.\ the argument for \eqref{convergence-of-sum}) 
\[
\phi(\sigma) = \Phi_\sigma(\wt u_\sigma) 
= \lim_{k\to \infty} \phi^{\eta^\sigma_k}(\sigma)
= \lim_{k\to \infty} \Phi^{\eta^\sigma_k}_\sigma\big( \tyos
u\sigma{\eta^\sigma_k} \big).  
\]

\noindent
\STEP{5: The simple De  Giorgi estimate.} Here we use the trivial lower
bound $ \wrslopar{\eta}{w} {\rho} \geq \calS_{\calR_\eta} (w) \geq \rslo w$
because of $\calR^*_\eta\geq \calR^*$. We start from
\eqref{eq:DeGiIdent-eta-slope} 
and obtain for all fixed $\sigma\in (0,\sigma_*)$ and all $\eta\in (0,1)$ the
estimates 
\[
\Phi_\sigma^\eta\big( \tyos u \eta\sigma \big)
 + \int_0^\sigma \rslo { \tyos u\rho\eta } \dd \rho \leq \calE(\ugi).
\]
In Step 4 we have shown that the first term converges to the desired
limit. For the integral it suffices to use Fatou's lemma and the 
liminf estimate \eqref{eq:SR.liminf}. Thus, we obtain
\[
\Phi_\sigma\big(\wt u_\sigma \big) 
 + \int_0^\sigma \rslo { \wt u_\rho } \dd \rho \leq \calE(\ugi),
\]
which is the desired simple De Giorgi estimate. 
\smallskip 

\noindent 
\STEP{6: The improved De Giorgi estimate.} We proceed as in the previous
steps, but apply the selection theorem to the enlarged family
$(g^\eta)_{\eta\in (0,1]}$ with  
\[
g^\eta:(0,\sigma_*) \to M_2  \quad \text{with } \ 
g^\eta_\sigma:=(\tyos u\sigma\eta, \eta) \in M_2:=  (B_{\ol E}(0){\cap} S_{\ol E} )  \ti
[0,1] \subset X\ti \R. 
\]
The measurable function $\calA$ is defined via the conditioned
$\calR_\eta$ slope as
\[
\calA\big(\sigma, (u,\eta) \big) := \left\{\begin{array}{cl} 
 \wrslopar{\eta}u\sigma  & \text{for } u\in \Jugi_\eta(\sigma), \\ 
\infty& \text{otherwise},   \end{array}
\right. 
\]
where for $\eta=0$ we use the notations $\calR_0=\calR$ and
$\Jugi_0(\sigma) = \Jugi(\sigma)$.  The additional Hypothesis
\eqref{eq:LSC.CR} is now exactly tailored in such a way that $\calA$ is lower
semi-continuous on $(0,\sigma_*)\to M_2$. Hence, Proposition
\ref{pr:MeasSelect} provides a selection $\wt u$ such that
\[
 \tyos u\sigma\eta \weak \wt u_\sigma \ \text{and } \  
  \wrslo {\wt u_\sigma}\sigma \leq  
\liminf_{\eta\to 0}  \wrslopar{\eta}{\tyos u\sigma\eta}\sigma \quad \text{for
  all } \sigma\in (0,\sigma_*).
\]
As in Step 4 and 5 we can pass to the liminf in \eqref{eq:DeGiIdent-eta-slope}
by combining the last estimate and Fatou's lemma, and the improved De Giorgi
estimate follows.

 Finally, the improved estimate with selection \ref{eq:BGS.DeGiIneq-wsel}
directly follows via Lemma \ref{l:conseq-of-closedness}.
\end{proof}

\subsection{An alternative route to equality via the chain rule}
\label{ss:CR} 

In this section, we extend an argument sketched in the introduction by showing 
that, even without radial differentiability for the dissipation
potential $\calR$,  the improved De Giorgi estimate with force selection
\eqref{eq:BGS.DeGiIneq-wsel}  can be enhanced  to an identity along a
curve $\sigma \mapsto \wt u_\sigma$
\begin{enumerate}
\itemsep-0.2em
\item with suitable  regularity 
\item such that the chain rule holds for $\sigma \mapsto \calE(\wt u_\sigma)$.
\end{enumerate}
Thus, the motivation for this section is to support our conjecture 
that  the  De Giorgi identity is valid in more general situations 
than those understood by now.

We substantiate the above requirements in the following hypothesis and
highlight that, for the variational interpolant it is in fact sufficient to
have (piecewise) absolute continuity, whereas the results from Section
\ref{ss:regularity-added} even granted, under additional assumptions, (local)
Lipschitz continuity of the curve $\sigma \mapsto \wt u_\sigma$.

\begin{hypothesis}[Regularity of  the   variational interpolant\,\&\,chain rule]
\label{h:CR}
\leavevmode\newline
Let $\sigma_*>0$.  We consider a curve
$(0,\sigma_*) \ni \sigma \mapsto \wt u_\sigma\in J_\sigma (\ugi)$ such that
\begin{enumerate}
\item $\wt u$ is \emph{piecewise absolutely continuous} on $(0,\sigma_*)$,
  namely
  \begin{equation}
    \label{pcw-AC}
    \begin{aligned}
      \hspace{-0,8cm} \exists\, p \geq 1 \ \exists\, \text{a partition } \{
      \uptau_j \}_{j=0}^J \text{ of } [0,\sigma_*) \ \forall \, j \in \{1,..,
      J\}\colon \ \wt u \in \mathrm{AC}^p \big((\uptau_{j-1}, \uptau_j) ;
      \Spx\big)\,.
    \end{aligned}
  \end{equation}
  As a consequence, for all $ \{1,\ldots, J\}$ the one-sided limits
  $ \llim{\wt u_{\uptau_j}}:= \lim_{\sigma \to \llim{\uptau_j}} \wt u_\sigma$
  and $ \rlim{\wt u_{\uptau_{j-1}}}:= \lim_{\sigma \to \rlim{\uptau_{j-1}}} \wt
  u_\sigma$ exist;
\item there exists a measurable selection $\wt \xi : (0,\sigma_*) \to \Spx^*$
  with $\wt \xi_\sigma \in \pl \calE (\wt u_\sigma)$ for a.a.\
  $\sigma \in (0,\sigma_*) $ such that
  \begin{equation}
    \label{chain-rule}
    \begin{aligned}
      \forall\, j \in \{1, \ldots, J\} \colon \ & \sigma \mapsto \calE(\wt
      u_\sigma) \text{ is absolutely continuous on } ( \uptau_{j-1}, \uptau_j),
      \\
      & \text{and } \frac{\dd}{\dd \sigma} \calE(\wt u_\sigma) 
       = \langle \wt \xi_\sigma, \wt u_\sigma' \rangle \quad \foraa\; \sigma \in (\uptau_{j-1},
      \uptau_j) \,.
    \end{aligned}
  \end{equation}
\end{enumerate}
\end{hypothesis}

 Let us dwell on the chain-rule condition in Hypothesis \ref{h:CR}. First
of all, we mention that, in the proof of Theorem \ref{thm:est2id} we shall in
fact apply \eqref{chain-rule} to a measurable force selection
$\sigma \mapsto \wt\xi_\sigma\in \EL {\sigma}{\wt u_\sigma} $.   We also
emphasize the chain rule for $\calE$ \emph{evaluated along} the (assumedly)
absolutely continuous curve $\sigma \mapsto \wt u_\sigma$ is required, only.
Nonetheless, it is natural to wonder for which classes of energies the chain
rule holds \emph{in general}.  Some sufficient conditions for its validity,
among which $\lambda$-convexity of $\calE$ for some $\lambda \in \R$, were
provided in \cite[Prop.\,2.4]{MiRoSa13NADN}, \cite[Prop.\,A.1]{MiRo23}.  There,
it was shown that for \emph{any} pair
$(\mathsf{u}, \serifxi) \in \AC([a,b];\Spx) {\times} \rmL^1 ([a,b];\Spx^*)$
 satisfying 
\begin{equation}
\label{conditions-for-chain-rule}
\sup_{s\in [a,b]}\calE(\mathsf{u}(s))<\infty \quad \text{and} \quad
\int_{a}^{b} \|\serifxi(s)\|_{*} \, \| \mathsf{u}'(s)\| \dd s <\infty\,,
\end{equation}
 the energy  $s \mapsto \calE(\mathsf{u}(s))$ is absolutely
continuous and the chain rule formula
$\frac{\dd }{\dd s} (\calE {\circ} \mathsf{u}) = \langle \serifxi, \mathsf{u}'
\rangle $ holds a.e.\ in $(a,b)$.  For instance, the above chain rule holds for
the energy functional from Example \ref{PDE-example}, which is indeed
$\lambda$-convex for $\lambda =-C_W$, cf.\
\cite[Sec.\,4.2]{Mielke-Rossi-Stephan}.
\par Now, in the upcoming Theorem \ref{thm:est2id} we will consider curves
$\sigma\mapsto \wt u_\sigma\in J_\sigma (\ugi)$  and force selections
$\sigma \mapsto \wt{\xi}_\sigma\in \EL{\sigma}{\wt u_\sigma} $  as in
Hypothesis \ref{h:CR} \emph{and}  satisfying the  improved De Giorgi
estimate \eqref{eq:BGS.DeGiIneq-wsel} (see Theorem \ref{th:DGL.NON.RD.Simple}
for the existence of such selections).   Then, we will clearly have the
energy bound $\sup_{\sigma \in ]0,\sigma_*)} \calE (\wt u_\sigma) <\infty $, as
well as
\[
\int_0^{\sigma_*}  \calR^*({-}\wt\xi_\rho)  \dd \rho <\infty\,.
\]
Therefore, \emph{if} the dual dissipation potential $\xi \mapsto \calR^*(\xi)$
controls $ \| \xi \|_{*}^{p'}$ with $p'$ conjugate to $p$, then combining the
above estimate with the condition
$u \in \mathrm{AC}^p \big((\uptau_j, \uptau_{j+1}) ; \Spx \big)$, we conclude
that for the pair $(\tilde u,\tilde \xi)$  also the second estimate
in \eqref{conditions-for-chain-rule} holds.

We are now in a position to state the our last result that  the  De
Giorgi identity holds even without radial differentiability, if we have some
regularity of the variational interpolant
$\sigma \mapsto (\wt u_\sigma, \wt\xi_\sigma)$.

\begin{theorem}[From  the  De Giorgi estimate to  the identity] 
\label{thm:est2id} \mbox{} 
Let the $\GBGS$ \linebreak[4] $(X,\calE,\calR)$ satisfy Hypotheses
\ref{h:X}.  and let $\ugi\in \mafo{dom}(\calE)$ be fixed.  Let
$(0,\sigma_*) \ni \sigma \mapsto \wt u_\sigma\in J_\sigma $ and
$(0,\sigma_*) \ni \sigma \mapsto \wt{\xi}_\sigma\in \EL {\sigma}{\wt u_\sigma}$
 satisfy Hypothesis \ref{h:CR}, as well as the improved De Giorgi
estimate with selection \eqref{eq:BGS.DeGiIneq-wsel} on $(0,\sigma_*)$.

Then,  \eqref{eq:BGS.DeGiIneq-wsel}    holds as an equality, i.e.\
\begin{equation}
  \label{DG-eq-sigma}
  \calE(\wt u_{\sigma_*}) + {\sigma_*}\calR\big( \frac1{\sigma_*}(\wt
  u_{\sigma_*}{-}\ugi)\big) + \int_0^{\sigma_*} \calR^*(-\wt\xi_\rho) \dd \rho =
  \calE(\ugi)
\end{equation}
and, a fortiori, we have
\begin{equation}
  \label{minimality-tildexi}
  \wt \xi_\sigma \in  \mathfrak{A}_{\calR} (\sigma, \wt{u}_\sigma)=
  \mafo{argmin} \bigset{\calR^*({-}\xi)}{\xi\in \EL\sigma{\wt{u}_\sigma}}
  \ \  \foraa\, \sigma \in (0,\sigma_*)\,,
\end{equation}
hence  the  De Giorgi identity \eqref{MEI-Ban-cond} holds.
\end{theorem}
\begin{proof} 
\STEP{1:} Let us fix $j \in \{1,\ldots  , J{-}1\}  $. It is
immediate to check that $\lrlim{\wt u_{\uptau_j}} \in J_{\uptau_j} $, and
thus $\Phi_{\uptau_j}( \llim{\wt u_{\uptau_j}} )  =\phi( \uptau_j)
 = \Phi_{\uptau_j}( \rlim{\wt u_{\uptau_j}} ) $, i.e.\
\begin{equation}
\label{2blater}
\calE(\llim{\wt u_{\uptau_j}}) + \uptau_j \calR\Big( \frac1{\uptau_j}
(\llim{\wt u_{\uptau_j}}{-}\ugi)\Big) =   
\calE(\rlim{\wt u_{\uptau_j}}) + \uptau_j \calR\Big( \frac1{\uptau_j}
(\rlim{\wt u_{\uptau_j}}{-}\ugi)\Big) \,. 
\end{equation}

Let now $\sigma \mapsto \wt{\xi}_\sigma\in \EL {\sigma}{\wt u_\sigma}$ be
as in the statement.  Hypothesis \ref{h:CR} enables us to apply the chain rule
on any interval $[s_*, s^*] \subset (\uptau_j,\uptau_{j+1})$, thus concluding
that
\[
\lim_{\sigma \to \uptau_{j+1}^-}\calE( \wt u_{\sigma}) -\lim_{\sigma \to \uptau_{j}^+}\calE( \wt u_{\sigma})   = \int_{\uptau_{j}}^{\uptau_{j+1}} \langle \wt \xi_\rho, \wt u_\rho' \rangle \dd \rho\,.
\]
Now, we have that 
\[
  \lim_{\sigma \to \uptau_{j}^+}\calE( \wt u_{\sigma}) = \lim_{\sigma \to
    \uptau_{j}^+} \!\!\left( \phi(\sigma) {-} \sigma \calR \Big(
    \frac1{\sigma}(\wt u_\sigma {-} \ugi)\Big)\right) = \phi(\uptau_j) {-}
  \uptau_j \calR \Big( \frac1{\uptau_j}(\wt u_{\uptau_j} {-} \ugi) \Big) =
  \calE(\rlim{\wt u_{\uptau_{j}}})
\]
by the continuity of $\phi$ (recall Proposition \ref{prop:DeriMarg}) and of
$\calR$. Analogously,
$\lim_{\sigma \to \uptau_{j+1}^-}\calE( \wt u_{\sigma}) = \calE(\llim{\wt u_{\uptau_{j+1}}})$. Therefore,
\begin{equation}
\label{CR-e}
\calE(\llim{\wt u_{\uptau_{j+1}}}) -  \calE(\rlim{\wt u_{\uptau_{j}}}) =
\int_{\uptau_{j}}^{\uptau_{j+1}} \langle \wt \xi_\rho, \wt u_\rho' \rangle \dd
\rho\,. 
\end{equation}
 In turn, since $ \wt{\xi}_\sigma\in \EL {\sigma}{\wt u_\sigma}$  implies
$-\wt \xi_\sigma \in \partial \calR \left(\frac1\sigma (\wt u_\sigma{-}\ugi)
\right) $, we may apply the  chain rule for the functional
$\sigma \mapsto \sigma \calR \left(\frac1\sigma (\wt u_\sigma{-}\ugi) \right) $ and   find that
\begin{equation}
\label{CR-R}
\begin{aligned}
&
\uptau_{j+1} \calR\Big( \frac1{\uptau_{j+1}} (\llim{\wt u_{\uptau_{j+1}}}{-}\ugi)\Big) 
 - \uptau_j \calR\Big( \frac1{\uptau_j} (\rlim{\wt u_{\uptau_j}}{-}\ugi)\Big)
\\
& = \int_{\uptau_{j}}^{\uptau_{j+1}}  
 \Big(
  \calR \big(\frac1\rho (\wt u_\rho{-}\ugi) \big) + 
  \rho\, \big\langle {-}\wt \xi_\rho, \frac1\rho \wt u_\rho' {-} \frac1{\rho^2}
  (\wt u_\rho{-}\ugi) \big\rangle  \Big) \dd \rho\,.
  \end{aligned}
\end{equation} 
Adding \eqref{CR-e} and \eqref{CR-R}, observing the cancellation of the term
$ \int_{\uptau_{j}}^{\uptau_{j+1}} \langle \wt \xi_\rho, \wt u_\rho' \rangle
\dd \rho$ and rearranging the remaining integral terms, we find
\[
  \begin{aligned}
    & \calE(\llim{\wt u_{\uptau_{j+1}}}) + \uptau_{j+1} \calR\Big(
    \frac1{\uptau_{j+1}} (\llim{\wt u_{\uptau_{j+1}}}{-}\ugi)\Big) +
    \int_{\uptau_{j}}^{\uptau_{j+1}} \!\!  \Big( \big\langle {-}\wt \xi_\rho,
    \tfrac1\rho (\wt u_\rho{-}\ugi) \big\rangle - \calR \big(\frac1\rho (\wt
    u_\rho{-}\ugi) \big) \Big) \dd \rho
    \\
    & = \calE(\rlim{\wt u_{\uptau_{j}}}) + \uptau_j \calR\Big( \frac1{\uptau_j}
    (\rlim{\wt u_{\uptau_j}}{-}\ugi)\Big)\,.
  \end{aligned}
\]
Now, since $-\wt\xi_\rho \in \pl\calR\big(\frac1\rho(\wt u_\rho{-}\ugi) \big)$,
the integrand in the third term on the left-hand side equals
$ \calR^*({-}\wt \xi_\rho)$.  In turn, by \eqref{2blater}, the right-hand side
equals $\Phi_{\uptau_j}(\llim{\wt u_{\uptau_j}} )$. All in all, we
conclude that
\begin{equation}
\label{conclusion-for-j}
\begin{aligned}
& 
\calE(\llim{\wt u_{\uptau_{j+1}}})  + \uptau_{j+1} \calR\Big( \frac1{\uptau_{j+1}} 
 (\llim{\wt u_{\uptau_j}}{-}\ugi) \Big)  
+ \int_{\uptau_{j}}^{\uptau_{j+1}}  \calR^*({-}\wt \xi_\rho) \dd \rho
\\
& \qquad  = 
\calE(\llim{\wt u_{\uptau_{j}}})  + \uptau_{j} \calR\Big(
\frac1{\uptau_{j}} (\llim{\wt u_{\uptau_{j}}}{-}\ugi)\Big) \quad \text{for all
} j \in \{1,\ldots,  J{-}1\}  \,. 
\end{aligned}
\end{equation}
\medskip

\noindent
\STEP{2:} It remains to derive the analogue of \eqref{conclusion-for-j} for
$j=0$. With this aim, we consider the first interval $[0,\uptau_1]$ of the
partition and fix $\mu \in (0,\uptau_1]$. Repeating the arguments from Step 1
we find
\[
\begin{aligned}
& 
\calE(\llim{\wt u_{\uptau_{1}}})  + \uptau_{1} \calR\Big( \frac1{\uptau_{1}} 
(\llim{\wt u_{\uptau_{1}}}{-}\ugi) \Big)  
+ \int_{\mu}^{\uptau_{j+1}}  \calR^*({-}\wt \xi_\rho) \dd \rho
\\
&  = 
\calE(\wt u_{\mu})  + \mu \calR\Big( \frac1{\mu} (\wt u_{\mu}{-}\ugi)\Big) 
 = \Phi_\mu(\wt u_\mu) = \phi( \mu)  \to \calE( \ugi)\,,
\end{aligned}
\]
 where the last identity follows because $\wt u_\mu $ is a minimizer and
$\phi$ is the value function, and the convergence stems from
\eqref{added-last-mom}.  We thus conclude
\begin{equation}
  \label{conclusion-for-0}
  \calE(\llim{\wt u_{\uptau_{1}}})  + \uptau_{1} \calR\Big( \frac1{\uptau_{1}} 
  (\llim{\wt u_{\uptau_{1}}}{-}\ugi) \Big)  + \int_{0}^{\uptau_{j+1}}  
  \calR^*({-}\wt \xi_\rho) \dd \rho  = \calE(\ugi)\,.
\end{equation}
\medskip

\noindent
\STEP{3: Conclusion of the proof.} Adding \eqref{conclusion-for-j} for
$j=1, \ldots,  J{-}1$ and \eqref{conclusion-for-0} we obtain
\[
  \calE(\llim{\wt{u}_{\sigma_*}}) +\sigma_* \calR\Big( \frac1{\sigma_*}
  (\llim{\wt{u}_{\sigma_*}}{-}\ugi)\Big) + \int_{0}^{\sigma_*} \calR^*({-}\wt
  \xi_\rho) \dd \rho = \calE(\ugi)\,.
\]
 In analogy to \eqref{2blater} for $\tau_J=\sigma_*$ we also have
$\Phi_{\sigma_*}( \llim{\wt u_{\tau_J}}) = \phi(\sigma_*)=
\Phi_{\sigma_*}(\wt u_{\sigma_*}) $, such that \eqref{DG-eq-sigma} is
established.
\end{proof}

\paragraph*{Acknowledgments.} 
The authors are grateful to an anonymous reviewer for several improvements and
many helpful remarks. 
The research of A.M. has been partially
supported by Deutsche Forschungsgemeinschaft (DFG) through the Berlin
Mathematics Research Center MATH+ (EXC-2046/1, project ID: 390685689)
subproject ``DistFell''. R.R. acknowledges the support of GNAMPA (INDAM), and of the MIUR-PRIN Grant 2020F3NCPX 
``Mathematics for industry 4.0 (Math4I4)''.

{\small 

\bibliographystyle{alpha_AMs}

\begin{thebibliography}{11}\itemsep0.1em

\bibitem[AGS05]{AmGiSa05GFMS}
L.~Ambrosio, N.~Gigli, and G.~Savar{\'e}, \emph{Gradient flows in metric spaces
  and in the space of probability measures}, first ed., Lectures in Mathematics
  ETH Z\"urich, Birkh\"auser Verlag, Basel, 2005.

\bibitem[Amb95]{Ambr95MM}
L.~Ambrosio: \emph{Minimizing movements}. Rend. Accad. Naz. Sci. XL Mem. Mat.
  Appl. (5) \textbf{19} (1995) 1773--1799.

\bibitem[AmT04]{Ambrosio-Tilli}
L.~Ambrosio and P.~Tilli, \emph{Topics on analysis in metric spaces}, Oxford
  Lecture Series in Mathematics and its Applications, vol.~25, Oxford
  University Press, Oxford, 2004.

\bibitem[Asp68]{Asplund68}
E.~Asplund: \emph{Fr\'{e}chet differentiability of convex functions}. Acta
  Math. \textbf{121} (1968) 31--47.

\bibitem[Att84]{Att84VCFO}
H.~Attouch, \emph{Variational convergence for functions and operators},
  Applicable Mathematics Series, Pitman (Advanced Publishing Program), Boston,
  MA, 1984.

\bibitem[AuF90]{Aubin-Frankowska}
J.-P.~Aubin and H.~Frankowska, \emph{Set-valued analysis}, Birkh\"{a}user 1990.

\bibitem[Bac21]{Bacho:PHD}
A.~Bacho, \emph{On the nonsmooth analysis of nonlinear evolution inclusions of
  first and second order with applications}, Ph.D. thesis, Technischen
  Universit\"at Berlin, 2021.

\bibitem[BaP12]{BaPre12}
V.~Barbu and T.~Precupanu, \emph{Convexity and optimization in {B}anach
  spaces}, fourth ed., Springer Monographs in Mathematics, Springer, Dordrecht,
  2012.

\bibitem[Bre73]{Brez73}
H.\  Br\'ezis, \emph{Op\'erateurs maximaux monotones et semi-groupes de contractions dans les espaces de Hilbert}, 
  North-Holland Mathematics Studies, No.\ 5,  North-Holland Publishing Co., Amsterdam-London; American Elsevier Publishing Co., Inc., New York, 1973.

\bibitem[BBI01]{BuragoBuragoIvanov}
D.~Burago, Y.~Burago, and S.~Ivanov, \emph{A course in metric geometry},
  Graduate Studies in Mathematics, vol.~33, American Mathematical Society,
  Providence, RI, 2001.
  
  
\bibitem[CaV77]{Castaing-Valadier77}
C.~Castaing and M.~Valadier, \emph{Convex analysis and measurable
  multifunctions}, Lectures Notes in Mathematics, Vol. 580, Springer-Verlag,
  Berlin-New York, 1977.

\bibitem[De{\relax}93]{Degi93NPMM}
E.~De{\relax}~Giorgi, \emph{New problems on minimizing movements}, Boundary
  value problems for partial differential equations and applications, RMA Res.
  Notes Appl. Math., vol.~29, Masson, Paris, 1993, pp.~81--98.
  
  \bibitem[EkT99]{ET99}
  I.\ Ekeland and R.\ T\'emam, \emph{Convex Analysis and Variational Problems}, 
  Classics In Applied Mathematics, vol.~28,  Society for Industrial and Applied Mathematics (SIAM), Philadelphia,1999. 
  
\bibitem[IoT79]{IofTih79TEPe}
A.~D.~Ioffe and V.~M.~Tihomirov, \emph{Theory of extremal problems}, Studies in
  Mathematics and its Applications, vol.~6, North-Holland Publishing Co.,
  Amsterdam, 1979, Translated from the Russian by Karol Makowski.

\bibitem[Mie23]{Miel23IAGS}
A.~Mielke: \emph{An introduction to the analysis of gradients systems}.
  WIAS-Preprint 3022 (2023), arXiv:2306.05026.

\bibitem[MiR23]{MiRo23}
A.~Mielke and R.~Rossi: \emph{Balanced-viscosity solutions to
  infinite-dimensional multi-rate systems}. Arch. Ration. Mech. Anal.
  \textbf{247}:53 (2023) 1--100.

\bibitem[MRS13]{MiRoSa13NADN}
A.~Mielke, R.~Rossi, and G.~Savar{\'e}: \emph{Nonsmooth analysis of doubly
  nonlinear evolution equations}. Calc. Var. Partial Differential Equations
  \textbf{46}:1-2 (2013) 253--310.

\bibitem[MRS23]{Mielke-Rossi-Stephan}
A.~Mielke, R.~Rossi, and A.~Stephan: \emph{On time-splitting methods for
  gradient flows with two dissipation mechanisms}. Preprint arXiv:2307.16137
  (2023).
  
\bibitem[Mor06]{Mordu-I} B.~S.~Mordukhovich, \emph{Variational analysis and
    generalized differentiation. {I}},
  Grundlehren der mathematischen Wissenschaften [Fundamental Principles of
  Mathematical Sciences], Springer-Verlag, Berlin, 2006.
  
\bibitem[RMS08]{RMS08} R.~Rossi, A.~Mielke, and G.~Savar{\'e}: \emph{A metric
    approach to a class of doubly nonlinear evolution equations and
    applications}. Ann. Sc. Norm. Super.  Pisa Cl. Sci. (5) \textbf{7}:1 (2008)
  97--169.
  
\bibitem[RoSa06]{RossiSavare06} R.~Rossi and G.~Savar{\'e}: \emph{Gradient
    flows of non convex functionals in {H}ilbert spaces and
    applications}. ESAIM Control Optim. Calc. Var.  \textbf{12}:3 (2006)
  564--614.

\end{thebibliography}

\providecommand{\bysame}{\leavevmode\hbox to3em{\hrulefill}\thinspace}
\providecommand{\MR}{}

} 
\end{document}